\documentclass[12pt,a4paper]{article}
\usepackage{comment}
\usepackage{latexsym}
\usepackage{amsmath,amssymb,amsfonts}
\usepackage{epsfig,float}
\usepackage{authblk}
\usepackage{tikz}
\usepackage{marvosym}
\usepackage{ulem}   

\usepackage{graphicx}

\newtheorem{proposition}{Proposition}[section]
\newtheorem{theorem}[proposition]{Theorem}
\newtheorem{corollary}[proposition]{Corollary}

\newtheorem{definition}[proposition]{Definition}
\newtheorem{remark}[proposition]{Remark}

\newcommand{\R}{I\!\!R}

\newcounter{secnum}

\begin{document}

\title{Tracking disturbances in transmission networks} 

\author{J.-G. Caputo and A. Hamdi}

\maketitle

{\normalsize \noindent Laboratoire de Math\'ematiques, INSA de Rouen Normandie,\\
Normandie Universit\'e \\
76801 Saint-Etienne du Rouvray, France\\
E-mail: caputo@insa-rouen.fr, adel.hamdi@insa-rouen.fr\\
}

\date{}


\begin{abstract}
We study the nonlinear inverse source problem of
detecting, localizing and identifying
unknown accidental disturbances on
forced and damped transmission networks. A first result is
that strategic observation
sets are enough to guarantee detection of disturbances. To 
localize and identify them, we additionally need the observation
set to be absorbent. If this set is dominantly absorbent, then
detection, localization and identification can be done in
"quasi real-time". We illustrate these results with numerical
experiments. 
\end{abstract}



%
%

\section{Introduction}

Networks play a major role in the transmission of electricity and fluids such 
as gas, oil, or water. A canonical
example is the power grid, a remarkable engineering realization that emerged in the past century.
A natural mathematical model of these networks is a graph together with discrete
conservation laws. Kirchhoff's law and Ohm's law yield the discrete Laplacian: the DC load flow model in
the jargon of power grids. This is the model used by many grid operators.  Adding inertia to the system
one obtains the so-called graph wave equation. Adding driving and damping, we can model
the time evolution of an electrical grid \cite{Kuramoto,Machowski} and a fluid network \cite{Maas}
in a linearized approximation. This is a linearization of the so-called "swing equation" used to model
the electrical grid \cite{Pagnier2,Kou}.

The detection of disturbances on such networks is an important topic.
For the grid, disturbances are electromechanical oscillations
induced by faulty equipment. For other fluid network, they can be leaks, causing unexpected pressure losses.
Such faults can break equipment, see \cite{survey} for a survey 
of such events for the grid. Due to the importance of the problem,
many authors addressed the issue of identifying a faulty vertex in a network.
Nudell and Chakrabortty \cite{Nudell} propose a graph theoretic algorithm to find forced oscillation inputs; they
assume minimal information and estimate how the fault affects the Green's function of the Laplacian.
Chen et al \cite{chen13} use an energy-based method to locate oscillation sources in power systems.
Delabay et al \cite{Delabays23} propose to locate a fault using a maximum likelihood approach: here the
state of the network is assumed to be known everywhere, not the parameters.

In the literature, many works assume the unknown disturbances to 
have a given form and only parts of this (frequency, amplitude,...) 
are unknown and subject to the inverse problem.
To address these shortfalls, six years ago, we introduced a source 
detection approach \cite{chk18} for the graph wave
equation with no damping, assuming distinct eigenvalues of the graph 
Laplacian and a disturbance on a single vertex that vanishes 
before reaching the 
final monitoring time.
We defined a strategic observation set of vertices; this allows 
to reconstruct the state of the disturbance free system using only
measurements from this set. From these observations, one 
can locate the source using an overdetermined linear system 
based on an adjoint formulation. Interestingly, the rank of 
the matrix in this linear
system revealed graph conditions for placing observation vertices.
Finally the procedure permits to reconstruct 
the disturbance.

Here, we revisit the problem, allow the graph laplacian to 
admit eigenvalues of arbitrary multiplicities, include
damping in the dynamics and allow several disturbances. 
To monitor the network for disturbances and localizing-identifying
them, one can raise the following three questions
\begin{itemize}
\item How to select the observation set of vertices to detect the presence of a disturbance?
\item How to select the observation set to localize disturbances and
identify their intensities ?
\item Can this localization and identification be done in "real time " ?
\end{itemize}
We also need technical conditions for which these questions can be answered.

To tackle question one, we introduce the notion of {\it strategic observation
set} and show that observing such a set gives the state
of the network unequivocally. This leads to requiring condition 1 where 
we assume
a healthy time period wherein no disturbances occur. The detection
is done by examining the differences between the observations of the 
system and it's states when it is free of disturbances (healthy states).
Once disturbances have been detected, we use these residuals 
to localize and identify them.

Questions three is easier to address than question two. We
introduce the technical condition that disturbances occur
away from the observation set. Then we show that if
a submatrix of the graph Laplacian has full rank and the condition
above is met, then one can detect, localize and identify disturbances in 
"quasi-real time". We call such observation sets {\it dominantly absorbent},
for these, the number of vertices is larger than half the total 
number of network vertices. 
This is a strong requirement, in many cases only a limited number of 
sensors are available. We show that if the
observation set is {\it absorbent}, then detection, localization and identification
of disturbances remains possible but with a delay.
We present numerical experiments illustrating these results.\\
The article is organized as follows:
section 2 presents the mathematical model and states the inverse
problem,
section 3 shows how to detect accidental disturbances.
We present the inverse problem for dominantly 
absorbent observation sets and for absorbent observation sets in sections 4 
and 5 respectively. Section 6 illustrates these results with
numerical experiments and we conclude in section 7.

\section{Mathematical modeling and problem statement}
In the following, we consider the evolution of a miscible flow
on a network modeled by a graph ${\cal{G}}(\cal{V})$ where
$\cal{V}$ is a set of vertices. 
We assume there are $|{\cal V}|=N$ vertices and the graph is connected.

We consider the system of differential equations
\begin{eqnarray} 
\ddot X(t) + \eta \dot X(t) - \Delta X(t)=F(t),
\label{Linear_Eq}
\end{eqnarray}
where $X$ is the N dimensional vector defining the network state.
$\Delta$ is the $N \times N$ graph Laplacian matrix defined by
\begin{definition} \label{Laplacian}
The graph Laplacian matrix $\Delta$ is the real,
symmetric negative semi-definite matrix, such that \\
$\Delta_{kl} = \Delta_{lk} > 0$ ~~if $k$ and $l$ 
are connected, 0 otherwise, \\
$\Delta_{kk} =-\sum_{l\neq k} w_{kl},$ .
\end{definition}
Here, we choose $\Delta_{kl}=1$ for simplicity and to illustrate the
influence of the topology.

The parameter $\eta>0$ is a damping. 
Equation (\ref{Linear_Eq}) is established as a linearized model
of the swing-equations 
for the voltage on a power grid; there $x_n$ is associated to the
phase of the voltage at vertex $n$. Equation (\ref{Linear_Eq})
also models the evolution of the pressure in a fluid (water or gas) 
network in the linear limit \cite{cks13,Maas,gas}.

We now recall the notions of forward and inverse problems associated 
to (\ref{Linear_Eq}). \\
$\blacktriangleright$ {\bf{Forward problem:}}
Apart from disturbances, the right hand side term of equation 
(\ref{Linear_Eq}) is $F=F^{sour}$, a well known source vector that 
ensures a desired functioning of the network. Given a final monitoring 
time $T>0$, a source $F^{sour}_{n=1,\dots,N}\in L^2(0,T)$ and 
two initial network state conditions 
$X(0)=\big(x_1(0),\dots,x_N(0)\big)^{\top}$ and 
$\dot X(0)=\big(\dot x_1(0),\dots,\dot x_N(0)\big)^{\top}$, 
the problem (\ref{Linear_Eq}) admits a unique solution 
$X=\big(x_1,\dots,x_N\big)^{\top}$ whose $x_{n=1,\dots,N}\in H^2(0,T)$. 
This is the so-called {\it{forward problem}}. Thus, given a 
set ${\cal{S}}$ of network observation vertices, we 
define the following observation operator: 
\begin{eqnarray}
M[F]:=\big\{x_n(t) \;\; \mbox{in}\; (0,T),\; \forall n\in {\cal{S}}\big\}.
\label{observations}  
\end{eqnarray}

We therefore introduce the notions of observation and non-observation sets.
\begin{definition}\label{SandE}
The set ${\cal{S}}$ of observation vertices 
wherein the local network state is recorded is called observation set. The 
complementary of ${\cal{S}}$, ${\cal{E}}= \cal{V}-\cal{S}$ is 
the set of non-observation vertices. In the following we label
$N_{\cal{S}}= |{\cal{S}}|, ~~N_{\cal{E}}= |{\cal{E}}|$.
\end{definition}

$\blacktriangleright$ {\bf{Inverse problem:}} Here, we write
$$F(t)=F^{sour}(t) + F^{dis}(t), ~~\forall t \in (0,T)$$
where $F^{sour}$ is known and $F^{dis}$ is an unknown accidental
disturbance that disrupts the healthy state of the network. In the 
inverse problem, we are given time records $d_n(t)$ in $(0,T)$ of the local 
state $x_n(t)$ taken at observation vertices $n\in {\cal{S}}$. The 
task is to detect-localize-identify 
$F^{dis}$ from $d_n(t), n \in \cal{S}$ i.e. invert the map 
\begin{eqnarray}
M[F]=\big\{d_n(t) \;\; \mbox{in}\; (0,T),\; \forall n\in {\cal{S}}\big\}.
\label{aim}  
\end{eqnarray}

An important question while addressing this inverse problem is: 
Which observation set ${\cal{S}}$ ensures the identifiability of the 
unknown occurring disturbances $F^{dis}$? 
In the remainder, we show that $\cal{S}$ needs to be {\it{strategic}},
a spectral property. We will also require ${\cal{S}}$ to be 
absorbent i.e.
\begin{definition}
A set of vertices $\cal{S}$ of a graph is absorbent if each
vertex of the graph is directly connected to at least one vertex of $\cal{S}$.
\label{absorbent}
\end{definition}
or verify a stronger property by being {\it{dominantly absorbent}}
that will be introduced later. \\
See the appendix for an algorithm to find an absorbent
set $\cal{S}$ for a given graph.

\section{Detection of accidental disturbances}
To detect unknown disturbances that may affect the 
network, we request the observation set ${\cal{S}}$ to be {\it{strategic}}. 
This notion is based on the spectral properties of the graph Laplacian 
matrix $\Delta$ defined in (\ref{Laplacian}). Here, we introduce this 
notion and establish that for a healthy network i.e., $F^{dis}=0$, the 
local knowledge of its state in a {\it{strategic}} observation set 
determines in a unique manner its global state in all of its vertices. This 
assertion yields the main ingredient of our procedure to
detect disturbances.

\subsection{Strategic set of network vertices}
Since the notion of a {\it{strategic}} observation set of vertices 
depends on the spectral properties of the graph Laplacian matrix $\Delta$, 
we employ the following notations: 
\begin{itemize}
\item $K$ is the number of distinct eigenvalues of the matrix $\Delta$. 
Thus, $K\le N$.
\item For $k=1,\dots,K$: We denote by $-\omega_k^2$ the $k^{th}$ eigenvalue of the matrix $\Delta$, where the real numbers $\omega_{k=1,\dots,K}$ are such that $0=\omega_1<\omega_2< \dots < \omega_K$. Hence, the distinct eigenvalues of $\Delta$ are labeled in the following decreasing order:
\begin{eqnarray}
0=-\omega_1^2> -\omega_2^2 > \dots > -\omega_K^2.
\label{eigenValues}
\end{eqnarray}
\item $m_k$ is the multiplicity of $-\omega_k^2$ and $M=\displaystyle\max_{k\in \{1,\dots,K\}} m_k$. Then, $N=\sum_{k=1}^K m_k$.
\item For $k=1,\dots,K$: We use the indices $k_{\ell=1,\dots,m_k}$ to denote the normalized 
orthogonal eigenvectors $v^{k_{\ell}}$ associated to the $k^{th}$ eigenvalue $-\omega_k^2$ of the matrix $\Delta$:
\begin{eqnarray}
\Delta v^{k_{\ell}}=-\omega_k^2 v^{k_{\ell}}, \quad \mbox{for} \; \ell=1,\dots,m_k, \;\;  \mbox{where} \;k_1=k.
\label{eigenvectors}
\end{eqnarray}
\end{itemize}
Notice that if the graph is connected, there is a single zero 
eigenvalue i.e., $m_1=1$ \cite{crs01}. We now introduce 
the notion {\it{strategic}} set of vertices.
\begin{definition}
A set ${\cal{S}}$ of vertices is {\it{strategic}} if for 
all eigenvalues $-\omega^2_k$ of multiplicity $m_k$, there exists a 
subset $I_k$ of $m_k$ vertices in ${\cal{S}}$ defining 
a $m_k\times m_k$ invertible matrix:
\begin{eqnarray}
A^{(k)}_{i\ell}=v_i^{k_{\ell}}, \quad \mbox{for} \;\; \ell=1,\dots,m_k \quad \mbox{and}\quad \forall i\in I_k,
\label{Ak}
\end{eqnarray}
where $v^{k_{\ell}}\in\R^N$ are the $m_k$ eigenvectors 
associated to the eigenvalue $-\omega^2_k$ of $\Delta$.
\label{strategic_set} 
\end{definition}
Hence, ${\cal{S}}$ is {\it{strategic}} if for all eigenvalues $-\omega^2_k$ 
of multiplicity $m_k$, there exists a subset $I_k$ of $m_k$ row labels in 
${\cal{S}}$ such that the matrix $A^{(k)}$ whose columns are the 
$m_k$ eigenvectors $v^{k_{\ell}}$ is of full rank. \\
Note that a {\it{strategic}} set contains at least $M$ vertices, 
where $M$ is the maximal multiplicity of the eigenvalues.

This definition generalizes the one given in \cite{chk18} for the 
case where all eigenvalues of the Laplacian matrix $\Delta$ are 
of multiplicity $1$. It is also a transposition for discrete systems 
of the result in \cite{Jai}. \\
From Definition \ref{strategic_set}, it follows that every 
set containing a {\it{strategic}} set is also {\it{strategic}}. As 
examples: 
\begin{enumerate}
\item For $M=1: \; {\cal{S}}$ is {\it{strategic}} if: $\forall k\in \{1,\dots,K=N\}, \exists I_k=\{i\}\in {\cal{S}},  v_i^k \ne 0$.  
\item For $M=2: \; {\cal{S}}$ is {\it{strategic}} if: $\forall k\in \{1,\dots,K<N\}$, it holds:

\hspace{0.5cm} $\bullet$ If $m_k=1$, then there exists $I_k=\{i\}\in {\cal{S}}$ such that $v_i^k \ne 0$.

\hspace{0.5cm} $\bullet$ If $m_k=2$, there exists $I_k=\{i,j\}\subset {\cal{S}}$ defining
$A^{(k)}=\begin{pmatrix}
v_i^{k_1} & v_i^{k_2}\\
v_j^{k_1} & v_j^{k_2}
\end{pmatrix}$
invertible.
\end{enumerate}
In the light of Definition \ref{strategic_set}, we establish the 
following theorem on the main property of a {\it{strategic}} set of 
vertices. This property so-called {\it{state 
representative of healthy networks}} stipulates that, in the
absence of disturbances, the states of the vertices of a {\it{strategic}} set 
determine uniquely the global state of the network. This is the healthy
state of the network.

\begin{theorem} 
Let $X_0\in\R^N$ and $\bar X_0\in\R^N$ be unknowns and $T^0\in(0,T)$. Provided ${\cal{S}}$ is a {\it{strategic}} set of vertices, if the solution $X^{R}(t)=\big(x^{R}_1(t),\dots,x^{R}_N(t)\big)^{\top}$ of the problem:
\begin{eqnarray}
\left\{
\begin{array}{lll}
\ddot X^{R}(t) + \eta \dot X^{R}(t)  - \Delta X^{R}(t) =0 \quad \mbox{in} \; (0,T^0),\\
X^{R}(0)=X_0  \quad \mbox{and} \quad \dot X^{R}(0)=\bar X_0,
\end{array}
\right.
\label{Rweq_0T0}
\end{eqnarray}
fulfills $x^{R}_n(t)=0, \forall t\in (0,T^0), \forall n\in {\cal{S}}$, then the initial conditions $X_0=\bar X_0=0$ in $\R^N$.
\label{Th_InitialConditions}
\end{theorem}
{\bf{Proof.}} See the appendix.

Thus, assuming the monitored network remains disturbance free 
during a period of time $(0, T^0)$, an immediate consequence 
of Theorem \ref{Th_InitialConditions} is the following corollary:  
\begin{corollary} Let $T^0\in(0,T)$, ${\cal{S}}$ be {\it{strategic}} and $F^{sour}$ be known in $(0, T^0)$. The data
\begin{eqnarray}
d_n(t)=x_n(t), \; \forall t\in(0, T^0), \; \forall n\in {\cal{S}},
\label{data_initialCond}
\end{eqnarray}
determines uniquely the solution $X(t)=\big(x_1(t),\dots,x_N(t)\big)^{\top}$ of the problem:
\begin{eqnarray}
\left\{
\begin{array}{lll}
\ddot X(t) + \eta \dot X(t)  - \Delta X(t) =F^{sour}(t) \quad \mbox{in} \; (0,T^0),\\
X(0)=X_0  \quad \mbox{and} \quad \dot X(0)=\bar X_0,
\end{array}
\right.
\label{weq_0T0}
\end{eqnarray}
where the two initial conditions $X_0\in\R^N$ 
and $\bar X_0\in\R^N$ are unknown.
\label{Corollary}
\end{corollary}
{\bf{Proof.}} Given $F^{sour}$ in $(0, T^0)$, let $X^i=\big(x_1^i,\dots,x_N^i\big)^{\top}$ be the solution of the problem (\ref{weq_0T0}) with the initial conditions $X^i(0)=X_0^i$ and $\dot X^i(0)=\bar X_0^i$, for $i=1,2$. We denote $d_n^i(t)=x_n^i(t), \forall t\in(0, T^0), \forall n\in {\cal{S}}$ and $Y=\big(y_1,\dots,y_N\big)^{\top}$ defined by $Y=X^2-X^1$.

Then, assuming $d_n^2(t)=d_n^1(t), \forall t\in(0, T^0), 
\forall n\in {\cal{S}}$, it follows that the state $Y$ satisfies 
$y_n(t)=0, \forall t\in (0,T^0), \forall n\in {\cal{S}}$ and solves 
the problem (\ref{Rweq_0T0}) with the initial conditions 
$Y(0)=X_0^2-X_0^1$ and $\dot Y(0)=\bar X_0^2 - \bar X_0^1$. 
Theorem \ref{Th_InitialConditions} implies that $X_0^2=X_0^1$ 
and $\bar X_0^2=\bar X_0^1$. Hence, $Y=0$ in $(0, T^0)$ which 
means that $X^2=X^1$ in $(0, T^0)$. \hspace{1cm} $\blacksquare$

\subsection{Procedure for detecting disturbances}
In the remainder of this paper, we assume the observation set ${\cal{S}}$ 
to be {\it{strategic}} and the following 
condition regarding unknown disturbances $F^{dis}$ to hold
\begin{eqnarray}
\mathbf{\big(C1\big):}\quad \exists T^0\in(0,T) \;\;\mbox{such that} \; \; F^{dis}(t)=0, \; \forall t\in(0,T^0).
\label{H1}
\end{eqnarray}
Hence, according to Corollary \ref{Corollary}, the data $d_n(t)=x_n(t), \forall t\in(0, T^0), \forall n\in {\cal{S}}$ determines in a unique manner the unknown network initial conditions $X(0)\in \R^N$ and $\dot X(0)\in \R^N$ involved in (\ref{weq_0T0}). Afterwards, from solving the problem (\ref{weq_0T0}) with the already determined network initial conditions and $T$ instead of $T^0$, we compute the {\it{healthy}} network state $X^{H}=\big(x_1^H,\dots,x_N^H\big)^{\top}$ in $(0,T)$. Notice that $X^H$ represents the network state throughout the entire monitoring period of time $(0,T)$ if it would remain out of all disturbances.

Besides, we introduce the residual variable $X^R=X-X^H$ in $(0,T)$. From assuming the condition $\mathbf{\big(C1\big)}$ holds true, it follows that $X^R=\big(x_1^R,\dots,x_N^R\big)^{\top}$ is subject to:
\begin{eqnarray}
\left\{
\begin{array}{lll}
\ddot X^{R}(t) + \eta \dot X^{R}(t)  - \Delta X^{R}(t) =F^{dis}(t) \quad \mbox{in} \; \big(T^0, T\big),\\
X^{R}(T^0)=\dot X^{R}(T^0)=0.
\end{array}
\right.
\label{Rweq_T0T}
\end{eqnarray}
Since $x_n^R=x_n-x_n^H, \forall n\in{\cal{V}}$, the problem (\ref{Rweq_T0T}) leads in particular to
\begin{eqnarray}
\forall \bar T\in (T^0,T), \quad F^{dis}=0 \;\; \mbox{in} \; (T^0, \bar T) \; \implies\;   x_n=x_n^H \;\; \mbox{in} \; (T^0, \bar T), \forall n\in {\cal{S}}. 
\label{R_implication}
\end{eqnarray}
In view of (\ref{R_implication}), we summarize the main steps defining the procedure we developed to detect disturbances in the following algorithm:

\hspace{0cm}{\bf{\uwave{Algorithm1}: Detection of unknown disturbances.}}

{\bf{Assume:}} ${\cal{S}}$ a {\it{strategic}} set of vertices and condition $\mathbf{\big(C1\big)}$ holds.

{\bf{Data:}} $T^0\in(0,T)$, $F^{sour}$ in $(0,T)$, $d_{n\in {\cal{S}}}$ in $(0,T)$, $\Delta t>0$ and a tolerance $\varepsilon>0$.

\hspace{0.2cm}{\bf{\underline{Begin}}}

\hspace{0.2cm} {\bf{1.}} Use data $d_{n\in {\cal{S}}}$ in $(0,T^0)$ to determine $X(0)$ and $\dot X(0)$ of the problem (\ref{weq_0T0}).

\hspace{0.2cm} {\bf{2.}} Solve (\ref{weq_0T0}) with $T$ instead of $T^0$ to compute the {\it{healthy}} state $X^H$ in $(0,T)$.

\hspace{0.2cm} {\bf{3.}} For $\bar T=T^0$ to $T$ in steps of $\Delta t$.

\hspace{1cm} $\blacktriangleright$ Compute the residual: ${\cal{R}}=\displaystyle\frac{1}{N_{\cal{S}}} \displaystyle\sum_{n \in {\cal{S}}}\big\|d_n-x_n^H\big\|_{L^2(\bar T, \bar T +\Delta t)}$.

\hspace{1cm} $\blacktriangleright$ If ${\cal{R}}>\varepsilon$, then {\it{detected disturbances}} in $(\bar T , \bar T+\Delta t)$.  Break. 

{\bf{\underline{End}.}}

The statement (\ref{R_implication}) raises a question about the effectiveness of the detection procedure: If (\ref{R_implication}) is not an equivalence, then we might have $x_n=x^H_n, \forall n \in {\cal{S}}$ in $(T^0, \bar T)$ whereas $F^{dis}\ne 0$ in $(T^0, \bar T)$. To address this issue, note first that in practice detection is more useful when unknown disturbances $F^{dis}=\big(F^{dis}_1,\dots,F^{dis}_N\big)^{\top}$ affect inaccessible/non-observation network vertices $n\in{\cal{E}}={\cal{V}}-{\cal{S}}$. Thus, in the remainder of the article, we consider that
\begin{eqnarray}
\mathbf{\big(C2\big):}\quad F_n^{dis}=0 \;\; \mbox{in} \; (0,T),\; \forall n\in {\cal{S}}.
\label{H2}
\end{eqnarray}
Provided $\mathbf{\big(C2\big)}$ holds, selecting from the $N$ equations defining the problem (\ref{Rweq_T0T}) all those associated to labels of the observation vertices $n\in{\cal{S}}$, we obtain: For all $t\in (T^0,T)$,
\begin{eqnarray}
\begin{array}{llll}
\ddot x_n^R(t) + \eta \dot x_n^R(t) - \displaystyle\sum_{m\in{\cal{V}}} \Delta_{nm}x_m^R(t)=0,\; \forall n\in {\cal{S}} \quad \Leftrightarrow \quad \displaystyle\sum_{m\in{\cal{E}}} \Delta_{nm}x_m^R(t)=d^R_n(t),\; \forall n\in {\cal{S}}\\
\mbox{where} \;\; d^R_n(t)=\ddot x_n^R(t) + \eta \dot x_n^R(t) - \displaystyle\sum_{m\in{\cal{S}}} \Delta_{nm}x_m^R(t).
\end{array}
\label{Residual_LS}
\end{eqnarray}
Let $\Delta\big({\cal{S}};{\cal{E}}\big)$ be the $N_{\cal{S}}\times N_{\cal{E}}$ submatrix obtained by selecting from the $N \times N$ matrix $\Delta$ the $N_{\cal{S}}$ rows associated to the labels of the observations vertices $n\in{\cal{S}}$ and the $N_{\cal{E}}$ columns associated to the labels of the non-observation vertices $m\in{\cal{E}}$. Then, (\ref{Residual_LS}) reads
\begin{eqnarray}
\Delta\big({\cal{S}};{\cal{E}}\big)X^{R_{\cal{E}}}(t)=d^R(t),\; \forall t\in(T^0,T),
\label{Residual_LS_matrix}
\end{eqnarray}
where $d^R\in \R^{N_{\cal{S}}}$ is defined by $d_n^R, \forall n\in{\cal{S}}$ and $X^{R_{\cal{E}}}\in\R^{N_{\cal{E}}}$ is given by $x_n^R, \forall n\in {\cal{E}}$. The effectiveness of the detection procedure relies on the properties of the matrix $\Delta\big({\cal{S}};{\cal{E}}\big)$ involved in (\ref{Residual_LS_matrix}). In the remainder, depending on the structure of the network, we consider the following two types of strategic observation sets.  If the matrix $\Delta\big({\cal{S}};{\cal{E}}\big)$is full rank,  ${\cal{S}}$ is said to be {\it{dominantly absorbent}}. Otherwise, we will assume ${\cal{S}}$ to be absorbent.

\section{{\color{black}Inversion} using dominantly absorbent observations}

To ensure that the detection procedure is effective, we introduce the notion of {\it{dominantly absorbent}} set of vertices. Similarly to the {\it{strategic}} property that yields state representativity of {\it{healthy}} networks, the property {\it{dominantly absorbent}} ensures state representativity of {\it{contaminated}} networks in the sense that the assertion (\ref{R_implication}) holds as an equivalence. Thus, if $x_n(t)=x_n^H(t), \forall t \in (T^0, \bar T)$ for every vertex $n$ of a {\it{dominantly absorbent}} observation set, then $x_n(t)=x_n^H(t), \forall t \in (T^0, \bar T)$ for all vertices $n$ of the network. Hence, $F^{dis}=0 $ in $(T^0, \bar T)$ which means effectiveness of the detection procedure.

We start by introducing the notion of {\it{dominantly absorbent}}  set of vertices. Then, we prove that in addition to making the detection effective, this new notion also ensures identifiability and leads to  localize-identify all unknown detected disturbances.

\subsection{Definition and practical guidelines}
The notion of {\it{dominantly absorbent}} is defined as follows:
\begin{definition}
Let ${\cal{S}}$ be a subset of $N_{\cal{S}}$ vertices selected from the $N=N_{\cal{S}}+ N_{\cal{E}}$ vertices defining the set ${\cal{V}}={\cal{S}} + {\cal{E}}$ of all network vertices. ${\cal{S}}$ is {\it{dominantly absorbent}} if
the $N_{\cal{S}}\times N_{\cal{E}}$ matrix $\Delta\big({\cal{S}};{\cal{E}}\big)$ introduced in  (\ref{Residual_LS})-(\ref{Residual_LS_matrix}) is of full rank. 
\label{DominatlyAbsorbant} 
\end{definition}  
We point out the following remarks 
\begin{itemize}
\item The properties {\it{strategic}} and {\it{dominantly absorbent}} 
differ in terms of their definitions and purposes. \\ 
The notion of {\it{strategic}} depends on the spectral properties 
of the graph Laplacian matrix $\Delta$ whereas the notion of 
{\it{dominantly absorbent}} depends on the structure of the 
matrix $\Delta$ itself. \\
As for purpose, the {\it{strategic}} property ensures 
uniqueness of the unknown network initial conditions for a given 
source term in the healthy problem (\ref{weq_0T0}). The 
{\it{dominantly absorbent}} property ensures uniqueness of 
an unknown source contaminating the network state of 
homogeneous initial conditions in problem (\ref{Rweq_T0T}).
\item An immediate consequence of Definition \ref{DominatlyAbsorbant} is 
that 
$$N_{\cal{S}}\ge N_{\cal{E}} \;\Leftrightarrow\; N_{\cal{S}}\ge N/2 . $$
If this were not the case, 
$\Delta\big({\cal{S}};{\cal{E}}\big)$ would not have full rank.
\item A {\it Dominantly Absorbent} set is {\it Absorbent}. \\
To argue this assertion, 
note that if a subset $\cal{S}$ is {\it{absorbent}} 
(see Definition (\ref{absorbent}))  then, 
the graph Laplacian matrix $\Delta$ verifies
\begin{eqnarray}
\forall m\in {\cal{E}}, \; \exists n \in {\cal{S}} \quad \mbox{such that} \quad \Delta_{nm}\ne 0,
\label{Necessary_observ}
\end{eqnarray}
Then, we establish the following result:
\begin{proposition}
Let ${\cal{V}}$ be the set of all vertices defining a network whose graph Laplacian matrix is given in (\ref{Laplacian}). All {\it{dominantly absorbent}} subset ${\cal{S}} \subset {\cal{V}}$ is {\it{absorbent}}.
\label{Link_DA_A}
\end{proposition}
{\bf{Proof.}} We prove Proposition \ref{Link_DA_A} by contrapositive. Let ${\cal{S}}$ be a subset defined by at least half of the vertices in ${\cal{V}}$. Assume that ${\cal{S}}$ is not {\it{absorbent}}. Then, relation (\ref{Necessary_observ})
implies: 
\begin{eqnarray}
  \exists m_0\in {\cal{E}} \;\; \mbox{such that} \;\; \forall n \in {\cal{S}},\; \Delta_{nm_0}=0.
\label{Necessary_observ_consq}
\end{eqnarray}
The assertion (\ref{Necessary_observ_consq}) implies that the column $m_0$ of the matrix $\Delta\big({\cal{S}};{\cal{E}}\big)$ is null. Hence, $\Delta\big({\cal{S}};{\cal{E}}\big)$ cannot be of full rank which implies that it is not {\it{dominantly absorbent}}. \hspace{0cm}$\blacksquare$
\end{itemize}

In practice, selecting a {\it{dominantly absorbent}} observation set 
requires to choose $N_S \ge N/2$. This is a necessary condition
for the matrix $\Delta({\cal S};{\cal E})$ to be of full rank. However, it is
not sufficient. For this purpose, we recover a condition 
linked to the topology of the graph as in \cite{chk18}. Indeed, 
for a particular graph that has a {\it{joint}}, we establish  
a result that guides the selection of the observation vertices.

We recall what is a {\it{joint}} vertex in a graph:
\begin{definition}
In a graph, a joint is a vertex whose removal increases the number of 
connected components.
\end{definition}
Then on the selection of a {\it{dominantly absorbent}} set of vertices, we have:
\begin{theorem} Consider a graph admitting a joint at a vertex $k\ge N/2$ whose removal splits the graph vertices into two connected big and small 
sets. If the vertices defining ${\cal S}$ are all taken from the 
big set, then the matrix $\Delta\big({\cal{S}};{\cal{E}}\big)$ is 
not of full rank.
\label{jointtheorem}
\end{theorem}
The proof of Theorem \ref{jointtheorem} is very similar to the one established in our previous work \cite{chk18}.

\begin{figure} [H]
\centerline{
\epsfig{file=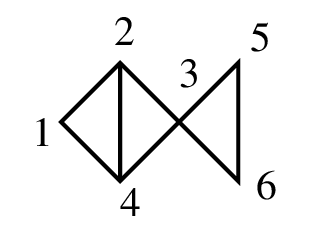,height=4 cm,width=12 cm,angle=0}
}
\caption{A graph with a joint.}
\label{joint}
\end{figure}

To illustrate this result, we consider the $6-$vertex graph 
admitting a {\it{joint}} at $k=3$ shown in Fig. \ref{joint}
The $6\times 6$ graph Laplacian matrix is given by:
\begin{eqnarray}
\Delta = 
\begin{pmatrix}
-2  & 1 & 0 & 1 & 0 & 0 \\
1   &-3 & 1 & 1 & 0 & 0 \\
0   & 1 &-4 & 1 & 1 & 1 \\
1   & 1 & 1 & -3& 0 & 0 \\
0   & 0 & 1 & 0 &-2 & 1 \\
0   & 0 & 1 & 0 &1  & -2 
\end{pmatrix}.
\end{eqnarray}
The choice
${\cal S}= \{1,2,3,4\}$ and ${\cal E}=\{5,6\}$ is such that
$N_S \ge N/2$, however the matrix
$$\Delta\big({\cal{S}};{\cal{E}}\big) =
\begin{pmatrix}
0 & 0 \\
0 & 0 \\
1 & 1 \\
0 & 0 
\end{pmatrix}$$
is not of full rank 2.
On the other hand,  ensuring that ${\cal{S}}$ contains at least one vertex from each graph connected component, for example, ${\cal{S}}= \{1,2,5\}$ and ${\cal{E}}=\{3,4,6\}$ leads to:
\begin{eqnarray}
\Delta\big({\cal{S}};{\cal{E}}\big) =
\begin{pmatrix}
0 &1 & 0 \\
1 &1 & 0 \\
1 &0 & 1 
\end{pmatrix}.
\label{Choice2}
\end{eqnarray}
We have $det(\Delta\big({\cal{S}};{\cal{E}}\big))=-1 $ so that the 
the set ${\cal{S}}= \{1,2,5\}$ is {\it{dominantly absorbent}}.

\subsection{Identifiability-Identification of detected disturbances}

In the light of Definition \ref{DominatlyAbsorbant}, we establish the following result on the identifiability of unknown disturbances $F^{dis}$ detected in $(T^0, \bar T)$ affecting the source $F=F^{sour}$ of the problem (\ref{Linear_Eq}) to become the undesired term $F=F^{sour}+F^{dis}$ in $(T^0,\bar T)$:
\begin{theorem}
Provided $F^{dis}=\big(F_1^{dis},\dots,F_N^{dis}\big)^{\top}$ 
fulfills ${\mathbf{\big(C1)}}$ and ${\mathbf{\big(C2)}}$, let 
$X=\big(x_1,\dots,x_N\big)^{\top}$ be the solution of the 
problem (\ref{Linear_Eq}) with $F=F^{sour}+F^{dis}$. If 
the {\it{strategic}} set ${\cal{S}}$ defining $M[F]$ in 
(\ref{observations}) is {\it{dominantly absorbent}}, then: 
the data
\begin{eqnarray}
d_n(t)=x_n(t), \; \forall t\in (0,T),\; \forall n\in{\cal{S}},
\label{data_T0T}
\end{eqnarray}
determines in a unique manner all detected unknown disturbances 
$F^{dis}_n\in L^2(T^0, T)$. 
\label{Th_Identifiability}
\end{theorem} 
{\bf{Proof.}} Let 
$X^{(i)}=\big(x_1^{(i)},\dots,x_N^{(i)}\big)^{\top}$ be the solution of the problem (\ref{Linear_Eq}) with $F^{(i)}=F^{sour}+F^{{dis}^{(i)}}$ in $(0, T)$, for $i=1,2$. Let $d^{(i)}_n(t)=x^{(i)}_n(t), \forall t\in (0,T), \forall n\in{\cal{S}}$ and $Y=\big(y_1,\dots,y_N\big)^{\top}$ be the state defined by $Y=X^{(2)} - X^{(1)}$. Assume that
\begin{eqnarray}
d^{(2)}_n(t)=d^{(1)}_n(t), \; \forall t\in (0,T),\; \forall n\in{\cal{S}}.
\label{data_T0T_ident}
\end{eqnarray}
Provided $F^{{dis}^{(i=1,2)}}$ fulfill ${\mathbf{\big(C1)}}$, the state $Y$ solves the problem (\ref{Rweq_0T0}) with the initial conditions $Y(0)=X^{(2)}(0)-X^{(1)}(0)$ and $\dot Y(0)=\dot X^{(2)}(0)- \dot X^{(1)}(0)$. Since the set ${\cal{S}}$ is {\it{strategic}} and (\ref{data_T0T_ident}) leads to $y_n(t)=0, \forall t\in(0, T^0), \forall n\in{\cal{S}}$, from applying Theorem \ref{Th_InitialConditions} it follows that $Y(T^0)=\dot Y(T^0)=0$. Therefore, the state $Y$ is subject to:
\begin{eqnarray}
\left\{
\begin{array}{lll}
\ddot Y(t) + \eta \dot Y(t)  - \Delta Y(t) =F^{{dis}^{(2)}}(t)-F^{{dis}^{(1)}}(t) \quad \mbox{in} \; \big(T^0, T\big),\\
Y(T^0)=\dot Y(T^0)=0.
\end{array}
\right.
\label{Yweq_T0T}
\end{eqnarray}
From (\ref{data_T0T_ident}), it follows that $y_n(t)=0, \; \forall t\in (T^0,T), \forall n\in{\cal{S}}$. Provided $F^{{dis}^{(i=1,2)}}$ fulfill ${\mathbf{\big(C2)}}$, by applying to the problem (\ref{Yweq_T0T}) similar techniques as in (\ref{Residual_LS})-(\ref{Residual_LS_matrix}) we get 
\begin{eqnarray}
\Delta\big({\cal{S}};{\cal{E}}\big)Y^{{\cal{E}}}(t)=0,\; \forall t\in(T^0, T),
\label{Identifiability_LS_matrix}
\end{eqnarray}
where $Y^{{\cal{E}}}\in\R^{N_{\cal{E}}}$ defined by $y_n, \forall n\in {\cal{E}}$. Since the set ${\cal{S}}$ is {\it{dominantly absorbent}}, the matrix $\Delta\big({\cal{S}};{\cal{E}}\big)$ is of full rank. Hence, (\ref{Identifiability_LS_matrix}) implies that $Y^{{\cal{E}}}=0$ in $(T^0,T)$ which leads to $Y=0$ in $(T^0,T)$. From the problem (\ref{Yweq_T0T}), it follows that $F^{{dis}^{(2)}}=F^{{dis}^{(1)}}$ in $(T^0, T)$. \hspace{0cm}$\blacksquare$ 

To identify the unknown disturbances $F^{dis}$ detected in $(T^0, \bar T)$, we recall the residual state $X^R=X-X^H$ that solves the problem (\ref{Rweq_T0T}). Since its part $X^{R_{{\cal{E}}}}$ containing the unknown local residuals $x_n^R, \forall n\in{\cal{E}}$ solves the linear system (\ref{Residual_LS_matrix}), we proceed as follows:

{\bf{\uwave{Algorithm2}:  Identification using a dominantly absorbent observation set.}}

{\bf{Assume:}}  ${\cal{S}}$ is {\it{strategic}} and {\it{dominantly absorbent}} and conditions $\mathbf{\big(C1\big)}$ and $\mathbf{\big(C2\big)}$ hold.

{\bf{Data:}} First detection instant $\bar T\in (T^0,T)$ and the observations $d_{n\in {\cal{S}}}$ in $(\bar T,T)$.

\hspace{0.2cm}{\bf{\underline{Begin}}}

\hspace{0.2cm} {\bf{1.}} From (\ref{Residual_LS}), use the data $d_{n\in {\cal{S}}}$ in $(\bar T, T)$ to compute $d_n^R, \forall n\in{\cal{S}}$ in $(\bar T, T)$.

\hspace{0.2cm} {\bf{2.}} Solve the linear system (\ref{Residual_LS_matrix}) to determine $X^{R_{{\cal{E}}}}$ in $(\bar T,T)$.

\hspace{0.2cm} {\bf{3.}} Using the determined $X^{R_{{\cal{E}}}}$ in $(\bar T,T)$, reconstruct from (\ref{Rweq_T0T}) the unknown $F^{dis}$.

{\bf{\underline{End}.}}

Algorithm2 leads to determine $F^{dis}$ almost in real time, depending on the sampling period $\Delta t$. Provided the invertible matrix extracted from $\Delta\big({\cal{S}};{\cal{E}}\big)$ is well conditioned, we expect  the reconstruction to be accurate.

\section{{\color{black}Inversion} using absorbent observations}

We address the inverse problem where the number of observation vertices is very limited smaller than half of the total network vertices. Hence, the {\it{strategic}} set ${\cal{S}}$ is not anymore {\it{dominantly absorbent}} in the sense of Definition \ref{DominatlyAbsorbant}. Consequently, the matrix $\Delta\big({\cal{S}};{\cal{E}}\big)$ defining the linear system (\ref{Residual_LS_matrix}) cann't be of full rank which implies no guarantee about the existence/uniqueness of its solution. This raises  a {\it{detection effectiveness issue}} since the assertion (\ref{R_implication}) might not hold as equivalence as well as an {\it{identifiability-identification issue}} because the solutions of (\ref{Identifiability_LS_matrix}) and (\ref{Residual_LS_matrix}) might be not what we expect.

In practice, disturbances usually bring the network out of its operational functioning state and always end up by affecting all of its vertices. Therefore, we believe the developed procedure to detect disturbances remains operational but it might deliver a delayed detection of the first instant disturbances. Moreover, failing the {\it{dominantly absorbent}} condition, we request the {\it{strategic}} set of observation vertices ${\cal{S}}$ to be at least {\it{absorbent}}.

In the remainder, we assume ${\mathbf{\big(C1)}}, {\mathbf{\big(C2)}}$ hold and the {\it{strategic}} set of vertices ${\cal{S}}$ to be not {\it{dominantly absorbent}} but rather {\it{absorbent}}. We develop an identification method to determine the unknown disturbances $F^{dis}$ detected at $\bar T\in(T^0,T)$. {\color{black}However, $\bar T$ obtained from Algorithm1 with ${\cal{S}}$ absorbent might be a delayed first instant of disturbances. Hence, we consider $\bar T_k\in (\bar T, T)$ and focus on identifying $F^{dis}$ in $(T^0, \bar T_k)$}. For this purpose, we need an orthogonal family of $L^2(0,\bar T_k)$ to expand the unknown states $x_m^R, \forall m \in {\cal E}$.

\subsection{Orthogonal family of Legendre polynomials}

We recall the Legendre orthogonal polynomials in $L^2(-1,1)$ defined by the following recurrence formula: For all $s\in(-1,1)$,
\begin{eqnarray}
\left\{
\begin{array}{llll}
P_{\ell+1}(s)=\displaystyle\frac{2\ell +1}{\ell +1} s P_{\ell}(s) - \displaystyle\frac{\ell}{\ell +1} P_{\ell-1}(s), \quad \forall \ell\ge 1,\\   
P_0(s)=1 \quad \mbox{and}  \quad P_1(s)=s.
\end{array}
\right.
\label{Legendre_1moins1}
\end{eqnarray}
Then, given $\bar T_k\in (\bar T, T)$ we employ the following change of variables: 
\begin{eqnarray}
t=\frac{\bar T_k}{2}\big(s + 1\big) \qquad \mbox{and} \qquad {\cal{P}}_{\ell}(t)=P_{\ell}(s)=P_{\ell}\big(\frac{2}{\bar T_k}t-1\big) ,
\label{changeVar_Legendre}
\end{eqnarray}
for all $s\in(-1,1)$ and $\ell\ge 1$. It follows that the sequence of Legendre's polynomials $\big\{{\cal{P}}_{\ell}\big\}_{\ell\ge 0}$ defined by (\ref{Legendre_1moins1})-(\ref{changeVar_Legendre}) forms an orthogonal basis of $L^2(0,\bar T_k)$, see \cite{Kwon,Liu}. For instance, the first four Legendre's polynomials are given by: For all $t\in(0,\bar T_k)$,
\begin{eqnarray}
\begin{array}{llll}
  {\cal{P}}_0(t)=1, \quad {\cal{P}}_1(t)=\frac{2}{\bar T_k} t - 1,\\
  {\cal{P}}_2(t)=\frac{3}{2}\big(\frac{2}{\bar T_k}t - 1\big)^2 - \frac{1}{2},\\
  {\cal{P}}_3(t)=\frac{5}{2}\big(\frac{2}{\bar T_k} t - 1\big)^3 - \frac{3}{2}\big(\frac{2}{\bar T_k} t - 1\big).
\end{array}
\label{Legendre_T0TD_L3}
\end{eqnarray}
We set the local states $x^R_m\in L^2(0,\bar T_k), \forall m\in{\cal{E}}$ of the residual variable $X^R=X-X^H$ to:
\begin{eqnarray}
x^R_m(t)=\displaystyle\sum_{\ell=0}^L \bar x_m^{\ell}{\cal{P}}_{\ell}(t), \;\; \forall t\in (0, \bar T_k),
\label{Expan_xRm}  
\end{eqnarray}
where $\bar x_m^{\ell}=\frac{1}{\big\|{\cal{P}}_{\ell}\big\|^2_{L^2(0, \bar T_k)}}\displaystyle\int_{0}^{\bar T_k}x^R_m(t){\cal{P}}_{\ell}(t)dt$ and $\|{\cal{P}}_{\ell}\|^2_{L^2(0,\bar T_k)}=\frac{\bar T_k}{2\ell +1}$, for $\ell=0,\dots,L$.

\subsection{Localization-Identification of detected disturbances}

We consider the integers $1<I_0<\bar I < I$ and the time step size $\Delta t=T/I$ to define the discrete times $t_i=i\Delta t$, for $i=1,\dots,I$ with $T^0=I_0\Delta t$ and $\bar T=\bar I \Delta t$. {\color{black} Then, we set $\bar T_k=\bar T + k \Delta t$, where $k \in \{1,\dots, I-\bar I\}$ with $\bar I_k = \bar I +k$ and focus on identifying $F^{dis}$ in $(T^0, \bar T_k)$}. From using (\ref{Expan_xRm}) in the linear system (\ref{Residual_LS}) taken at $t_{i=I_0+1,\dots,\bar I_k}$, we get:
\begin{eqnarray}
\displaystyle\sum_{m\in{\cal{E}}}\displaystyle\sum_{\ell=0}^L \Delta_{nm}{\cal{P}}_{\ell}(t_i) \bar x_m^{\ell}=d^R_n(t_i), \;\; \forall n\in{\cal{S}}, \;\; \mbox{for} \; i=I_0+1,\dots,\bar I_k.
\label{xRn_E}
\end{eqnarray}
For the sake of clarity, we label the $N_{{\cal{S}}}$ observation vertices of the set ${\cal{S}}=\big\{n_1,\dots,n_{N_{{\cal{S}}}}\big\}$ and the $N_{{\cal{E}}}$ non-observation vertices of the set ${\cal{E}}=\big\{m_1,\dots,m_{N_{{\cal{E}}}}\big\}$. Afterwards, the linear system (\ref{xRn_E}) can be written under the following matrix form:
{\small{
\begin{eqnarray}
{\cal{C}}{\cal{X}}^R={\cal{D}}^R, \qquad \mbox{where} \;\; {\cal{X}}^R=
\begin{pmatrix}
  \bar x_{m_1}^0\\
  \vdots\\
  \bar x_{m_1}^L\\
  \vdots\\
  \vdots\\
  \bar x_{m_{N_{\cal{E}}}}^0\\
  \vdots\\
  \bar x_{m_{N_{\cal{E}}}}^L\\
\end{pmatrix}\in\R^{(L+1)\times N_{\cal{E}}}
,\quad
{\cal{D}}^R=
\begin{pmatrix}
  d^R_{n_1}(t_{I_0+1})\\
  \vdots\\
  d^R_{n_1}(t_{\bar I_k})\\
  \vdots\\ 
  \vdots\\
  d^R_{n_{N_{\cal{S}}}}(t_{I_0+1})\\
  \vdots\\
  d^R_{n_{N_{\cal{S}}}}(t_{\bar I_k})\\
\end{pmatrix}\in\R^{\hat I\times N_{\cal{S}}}
\label{LinearSyst_YG}
\end{eqnarray}
}}
$\hat I=\bar I_k-I_0$ and the $\big(\hat I\times N_{\cal{S}}\big)\times\big((L+1)\times N_{\cal{E}}\big)$ matrix ${\cal{C}}$ is defined by
{\tiny{
\begin{eqnarray}
{\cal{C}}=
\begin{pmatrix}
\Delta_{n_1m_1}{\cal{P}}_0(t_{I_0+1}) & \dots & \Delta_{n_1m_1}{\cal{P}}_L(t_{I_0+1}) &\dots\dots&\Delta_{n_1m_{N_{\cal{E}}}}{\cal{P}}_0(t_{I_0+1}) & \dots & \Delta_{n_1m_{N_{\cal{E}}}}{\cal{P}}_L(t_{I_0+1})&\\
  \vdots & & \vdots & & \vdots & & \vdots&\\
\Delta_{n_1m_1}{\cal{P}}_0(t_{\bar I_k}) & \dots & \Delta_{n_1m_1}{\cal{P}}_L(t_{\bar I_k}) &\dots \dots&\Delta_{n_1m_{N_{\cal{E}}}}{\cal{P}}_0(t_{\bar I_k}) & \dots & \Delta_{n_1m_{N_{\cal{E}}}}{\cal{P}}_L(t_{\bar I_k})&\\
\vdots & & \vdots & & \vdots & & \vdots&\\
  \vdots & & \vdots & & \vdots & & \vdots&\\
  \Delta_{n_{N_{\cal{S}}}m_1}{\cal{P}}_0(t_{I_0+1}) & \dots & \Delta_{n_{N_{\cal{S}}}m_1}{\cal{P}}_L(t_{I_0+1}) &\dots\dots&\Delta_{n_{N_{\cal{S}}}m_{N_{\cal{E}}}}{\cal{P}}_0(t_{I_0+1}) & \dots & \Delta_{n_{N_{\cal{S}}}m_{N_{\cal{E}}}}{\cal{P}}_L(t_{I_0+1})& \\
\vdots & & \vdots & & \vdots & & \vdots&\\
\Delta_{n_{N_{\cal{S}}}m_1}{\cal{P}}_0(t_{\bar I_k}) & \dots & \Delta_{n_{N_{\cal{S}}}m_1}{\cal{P}}_L(t_{\bar I_k}) &\dots\dots&\Delta_{n_{N_{\cal{S}}}m_{N_{\cal{E}}}}{\cal{P}}_0(t_{\bar I_k}) & \dots & \Delta_{n_{N_{\cal{S}}}m_{N_{\cal{E}}}}{\cal{P}}_L(t_{\bar I_k})& 
\end{pmatrix}.
\label{C}
\end{eqnarray}
}}
According to condition ${\mathbf{\big(C1)}}$, the residual states $x_n^R=0, \forall n \in{\cal{V}}$ in $(0, T^0)$. We determine the unknown vector ${\cal{X}}^R$ subject to the linear system (\ref{LinearSyst_YG})-(\ref{C}) from solving:
\begin{eqnarray}
\displaystyle\min_{{\cal{X}}^R\in\R^{(L+1)\times N_{\cal{E}}}}\zeta\big({\cal{X}}^R\big):=\frac{1}{2}\Big\|{\cal{C}}{\cal{X}}^R - {\cal{D}}^R\Big\|_2^2 + \frac{\alpha}{2}\displaystyle\sum_{m=m_1}^{m_{N_{\cal{E}}}}\Big\|\displaystyle\sum_{\ell=0}^L \bar x_m^{\ell}{\cal{P}}_{\ell}\Big\|^2_{L^2(0,T^0)},
\label{Min_xi}
\end{eqnarray}
where $\|\cdot\|_2$ stands for the Euclidean norm and $\alpha>0$ is a regularization coefficient \cite{Donatelli}. 
{\color{black}
This coefficient has to be selected as a good compromise between 
taking into account the measurements during the time period of
disturbances $(T^0, \bar T_k)$ and 
requiring the states in the non-observation vertices to match 
their healthy references in $(0, T^0)$.}

Solving the minimization problem (\ref{Min_xi}) leads to determine an approximation in $(T^0, \bar T_k)$ of $X^R$ the solution of the residual problem (\ref{Rweq_T0T}). Then, by selecting from (\ref{Rweq_T0T}) all equations associated to the non-observation vertices, it follows that each detected unknown disturbance $F^{dis}_m$ occurring in a non-observation vertex $m\in{\cal{E}}$ is subject to: 
\begin{eqnarray}
F_m^{dis}(t)=\ddot x^R_m(t) + \eta\dot x^R_m(t) - \displaystyle\sum_{n=1}^N\Delta_{mn}x^R_n(t), \;\; \forall t\in(T^0, \bar T_k).
\label{mth_wave_eq}
\end{eqnarray}
Besides, for $i=I_0+1, \dots, \bar I_k$ we denote by $\varphi_i$ the impulse response solution of:
\begin{eqnarray}
\left\{
\begin{array}{lll}
\ddot \varphi_i(t) - \eta \dot \varphi_i(t)=\delta(t-t_i) \quad \mbox{in} \; \big(t_{i-1}, t_{i+1}\big),\\
\varphi_i(t_{i-1})=\varphi_i(t_{i+1})=0,
\end{array}
\right.
\label{ImpulseResponse_mi}
\end{eqnarray}
where $\delta(t-t_i)$ is the Dirac mass at $t_i$. Using the Heaviside function ${\cal{H}}$, it comes that
\begin{eqnarray}
  \varphi_i(t)=\displaystyle\frac{\displaystyle e^{\eta(t-t_i)} - 1}{\eta} {\cal{H}}\big(t-t_i\big) - \displaystyle\frac{e^{\eta\Delta t} - 1}{\eta\big(e^{2\eta\Delta t} - 1\big)}\Big(e^{\eta(t-t_{i-1})} - 1\Big), \;\;\forall t\in \big(t_{i-1}, t_{i+1}\big).
\label{varphi_i}
\end{eqnarray}
\begin{remark}
For $\eta \Delta t$ small enough, using $e^{x} - 1 \approx x$ it follows from (\ref{varphi_i}) that the impulse response $\varphi_i$ can be approximated by its affine version:
\begin{eqnarray}
\varphi_i(t)\approx
\left\{
\begin{array}{llll}
-\frac{1}{2}\big(t-t_{i-1}\big), \;\; \forall t\in \big(t_{i-1}, t_i\big),\\
\frac{1}{2}\big(t-t_{i+1}\big), \;\; \forall t\in \big( t_i, t_{i+1}\big).
\end{array}
\right.
\label{Affine_phi}
\end{eqnarray}
Graphically, (\ref{Affine_phi}) is the triangle of base $\big(t_{i-1}, t_{i+1}\big)$ and isosceles at $\varphi_i(t_i)=-\Delta t/2$.
\label{Remark_phi}
\end{remark} 
From multiplying (\ref{mth_wave_eq}) by $\varphi_i$ and integrating by parts on $\big(t_{i-1},t_{i+1}\big)$, we obtain 
\begin{eqnarray}
\displaystyle\int_{t_{i-1}}^{t_{i+1}}F_m^{dis}(t)\varphi_i(t)dt=x^R_m(t_i) - \Big[x^R_m(t)\dot\varphi_i(t)\Big]_{t_{i-1}}^{t_{i+1}} - \sum_{n=1}^N\Delta_{mn}\displaystyle\int_{t_{i-1}}^{t_{i+1}}x^R_n(t)\varphi_i(t)dt.
\label{Impulse_Reconst_Fm}
\end{eqnarray}
For the purpose of applying an appropriate numerical integration formula to determine the unknown disturbances $F_m^{dis}(t_i)$ from (\ref{Impulse_Reconst_Fm}), we establish the following result: 
\begin{proposition}
Let $a<b$ be two real numbers, $t_i=(a+b)/2$ and $\varphi_i$ be the function defined in (\ref{Affine_phi}) with $\big(a, b\big)$ instead of $\big(t_{i-1}, t_{i+1}\big)$. Then, 
\begin{eqnarray}
\displaystyle\int_a^b g(t)\varphi_i(t)dt=\displaystyle\frac{g(a) + 10g(t_i) + g(b)}{12} \varphi_i(t_i) \displaystyle\frac{b-a}{2},
\label{Own_integration}
\end{eqnarray}
\label{Own_numerical_scheme}
for all $g$ a real polynomial of degree smaller or equal to $2$ in $\big(a, b\big)$.
\end{proposition}
\hspace{0.2cm} {\bf{Proof.}} See the appendix.
\begin{remark}
To approximate the integrals in (\ref{Impulse_Reconst_Fm}), the numerical integration formula of Proposition \ref{Own_numerical_scheme} is more appropriate than Simpson's approximation that reads:
\begin{eqnarray}
\displaystyle\int_{t_{i-1}}^{t_{i+1}}F_m^{dis}(t)\varphi_i(t)dt \approx \displaystyle\frac{F_m^{dis}(t_{i-1})\varphi_i(t_{i-1}) + 4F_m^{dis}(t_{i})\varphi_i(t_{i}) + F_m^{dis}(t_{i+1})\varphi_i(t_{i+1})}{6} 2\Delta t.
\label{Simpson}
\end{eqnarray}
Since $\varphi_i(t_{i-1})=\varphi_i(t_{i+1})=0$, for constant $F_m^{dis}(t)=f_m^{dis}$ in $(t_{i-1}, t_{i+1})$ we get
\begin{eqnarray}
\begin{array}{lll}
\hspace{-2cm} \blacktriangleright \; \mbox{Proposition \ref{Own_numerical_scheme}}: \; \displaystyle\int_{t_{i-1}}^{t_{i+1}} F_m^{dis}(t)\varphi_i(t)dt=f^{dis}_m\varphi_i(t_i) \Delta t,\\
\hspace{-2cm} \blacktriangleright \; \mbox{Simpson's formula}: \; \displaystyle\int_{t_{i-1}}^{t_{i+1}} F_m^{dis}(t)\varphi_i(t)dt \approx \frac{4}{3}f^{dis}_m \varphi_i(t_i) \Delta t.
\end{array}
\label{Approximation}
\end{eqnarray}
Thus, Proposition \ref{Own_numerical_scheme} gives the exact value of the integral which is $f^{dis}_m$ times the area of the triangle defined by (\ref{Affine_phi}), whereas Simpson's approximation slightly overestimates it.
\end{remark}
From employing the numerical integration formula (\ref{Own_integration}) in the assertion (\ref{Impulse_Reconst_Fm}) and assuming ${\mathbf{\big(C1\big)}}$ holds, it follows that each detected unknown disturbance $F_m^{dis}$ is subject to:
\begin{eqnarray}
\left\{
\begin{array}{lll}
F_m^{dis}(t_{i-1}) + 10F_m^{dis}(t_i) + F_m^{dis}(t_{i+1})=\displaystyle\frac{12}{\Delta t \varphi_i(t_i)}d_m^i, \quad \mbox{for}\; i=I_0+1,\dots,\bar I_k,\\
F_m^{dis}(t_{I_0})=0,
\end{array}
\right.
\label{own_scheme}
\end{eqnarray}
where $t_{I_0}=T^0$ and the data $d_m^i$ is given by
\begin{eqnarray}
d_m^i=x^R_m(t_i) - \Big[x^R_m(t)\dot\varphi_i(t)\Big]_{t_{i-1}}^{t_{i+1}} -\displaystyle\frac{\Delta t \varphi_i(t_i)}{12}\displaystyle\sum_{n=1}^N\Delta_{mn}\Big(x^R_n(t_{i-1}) + 10x^R_n(t_i) + x^R_n(t_{i+1})\Big).
\label{data_own_scheme}
\end{eqnarray}
The benefit of (\ref{own_scheme})-(\ref{data_own_scheme}), as opposed to (\ref{mth_wave_eq}), is that they do not involve the derivatives of the residuals $x_n^R$. Since the residuals $x_m^R$ in the non-observation vertices $m\in{\cal{E}}$ are identified from the minimization problem (\ref{Min_xi}), the determination of their derivatives might add more computational errors and lead to misinterpretations.
\subsubsection{Localization of detected disturbances}

From (\ref{own_scheme})-(\ref{data_own_scheme}), by assuming $F_n^{dis}(t_j)=x_n^R(t_j)=0$, for $n=1,\dots,N$ and all $j\le i$ whereas $\exists m_0\in{\cal{E}}$ such that $F_{m_0}^{dis}(t_{i+1})\ne 0$ and $F_m^{dis}(t_{i+1})=0, \forall m\ne m_0$, it follows that
\begin{eqnarray}
\begin{array}{llll}
\bullet \; x_{m_0}^R(t_{i+1})=C_{m_0}\displaystyle\frac{\overbrace{F_{m_0}^{dis}(t_{i+1})}^{\mbox{\color{black}Active dis.}} + \overbrace{\displaystyle\sum_{n=1, n\ne m_0}^N \Delta_{m_0n}x_n^R(t_{i+1})}^{\mbox{\color{black}Passive dis.}}}{-\Delta_{m_0m_0}},\\
\bullet \; x_{m}^R(t_{i+1})=C_m\displaystyle\frac{\overbrace{\displaystyle\sum_{n=1, n\ne m}^N \Delta_{mn}x_n^R(t_{i+1})}^{\mbox{\color{black}Passive dis.}}}{-\Delta_{mm}},
\end{array}
\label{interpretation}
\end{eqnarray}
where $C_m=\Big(1 + \displaystyle\frac{12 \dot\varphi_i(t_{i+1})}{\Delta_{mm} \varphi_i(t_i)\Delta t}\Big)^{-1}, \forall m\in \{1,\dots,N\}$. The analysis of (\ref{interpretation}) leads to:
\begin{enumerate}
\item The residual $x_{m_0}^R$ in the vertex $m_0$ is the sum of both active and passive disturbances. The first is the actual disturbances $F_{m_0}^{dis}\ne 0$ affecting the vertex $m_0$ whereas the second is due to the coupling effects through the Laplacian matrix.
\item {\color{black}Since $\Delta_{mn}\in\{0,1\}, \forall n\ne m$ and $-\Delta_{mm}$ is the number of terms $\Delta_{mn}=1$, the residual $x_m^R$ in each vertex $m$ where $F_m^{dis}=0$ is defined only by the passive disturbance that is a normalized average of the residuals in all vertices directly connected to it.} 
\end{enumerate}
{\color{black}Therefore, each vertex $m_0\in{\cal{E}}$ hosting unknown actual disturbances is characterized by a residual $x_{m_0}^R$ that, as soon as $F_{m_0}^{dis}(t)\ne 0$, deviates from a common behavior followed by the residuals in all vertices $m$ where $F_m^{dis}=0$. That enables to localize $m_0$ in ${\cal{E}}$.}

\subsubsection{Identification of detected-localized disturbances}

According to (\ref{own_scheme})-(\ref{data_own_scheme}), from setting $F_m^{dis}(t_{\bar I_k +1})=F_m^{dis}(t_{\bar I_k})$, it follows that the unknown detected disturbances $F_m^{dis}(t_i)$, for $i=I_0+1,\dots, \bar I_k$ are subject to:
\begin{eqnarray}
\begin{pmatrix}
10 & 1  & 0  & 0 & 0 &\dots  & 0\\
1 & 10 & 1  & 0 & 0 & \dots & 0\\
&\vdots &  &  \vdots &  &\vdots &\\
0 & 0  & \dots & 0 & 1 & 10 & 1\\
0 & 0 &\dots&0 &0&1 & 11
\end{pmatrix}
\begin{pmatrix}
  F_m^{dis}(t_{I_0+1})\\
  \vdots\\
  \vdots\\
  \vdots\\
  F_m^{dis}(t_{\bar I_k})
\end{pmatrix}
=\frac{12}{\Delta t}
\begin{pmatrix}
d_m^{I_0+1}/\varphi_{I_0+1}(t_{I_0+1})\\
  \vdots\\
  \vdots\\
  \vdots\\
d_m^{\bar I_k}/\varphi_{\bar I_k}(t_{\bar I_k})
\end{pmatrix}.
\label{Lin_Syst_Fm}
\end{eqnarray}
From solving the linear system (\ref{Lin_Syst_Fm}) defined by a $\big(\bar I_k-I_0)\times(\bar I_k-I_0)$ tridiagonal Toeplitz matrix, we determine the sought disturbances $F_m^{dis}(t_i)$, for $i=I_0+1, \dots, \bar I_k$.

\hspace{-0.2cm}{\bf{\uwave{Algorithm3}:  Localization-Identification of detected disturbances}}

\hspace{0.2cm} {\bf{Assume:}}  ${\cal{S}}$ {\it{strategic}} and {\it{absorbent}} with conditions $\mathbf{\big(C1\big)}$ and $\mathbf{\big(C2\big)}$ hold.

\hspace{0.2cm} {\bf{Data:}} $T^0\in(0,T)$, $\bar T_k\in (\bar T, T)$ and the observations $d_{n\in {\cal{S}}}$ in $(0, \bar T_k)$.

\hspace{0.2cm}{\bf{\underline{Begin}}}

\hspace{0.2cm} {\bf{1.}} From (\ref{Residual_LS}), use $d_{n\in {\cal{S}}}$ in $(T^0, \bar T_k)$ to compute $d_n^R(t_i), \forall n\in{\cal{S}}$ for $i=I_0+1,\dots,\bar I_k$.

\hspace{0.2cm} {\bf{2.}} Solve the problem (\ref{Min_xi}) to determine $x_m^R(t_i), \forall m\in{\cal{E}}$ for $i=I_0+1,\dots,\bar I_k$.

\hspace{0.2cm} {\bf{3.}} {\bf{Localize}} vertices hosting disturbances by analyzing the behavior of $x_m^R, \forall m\in{\cal{E}}$. 

\hspace{0.2cm} {\bf{4.}}  {\bf{Identification:}} For each localized vertex $m\in{\cal{E}}$ hosting unknown disturbances:

\hspace{1.2cm} $\blacktriangleright$  From (\ref{data_own_scheme}), compute $d_m^i$ for $i=I_0+1,\dots,\bar I_k$.

\hspace{1.2cm} $\blacktriangleright$ Solve the linear system (\ref{Lin_Syst_Fm}) to determine $F_m^{dis}(t_{i})$, for $i=I_0+1,\dots,\bar I_k$.

\hspace{0.2cm}{\bf{\underline{End}.}}

\section{Numerical experiments}

\begin{figure} [H]
\centerline{
\epsfig{file=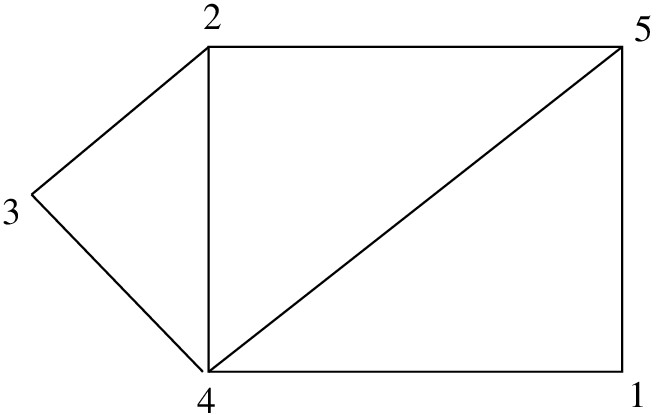,height=2 cm,width=8 cm,angle=0}
}
\caption{
A graph with five vertices.
}
\label{g5}
\end{figure}

We carry out numerical experiments on the $5-$vertex graph shown in Fig. \ref{g5}.
It's graph Laplacian matrix is given by:
$$\Delta = \left(
\begin{array}{ccccc}
-2 & 0 & 0 & 1 & 1 \\
0 & -3 & 1 & 1 & 1 \\
0 & 1  &-2 & 1 & 0 \\
1 & 1  & 1 & -4 & 1 \\
1 & 1  & 0 & 1  & -3 \\
\end{array}
\right). 
$$

\subsection{Determination of the healthy network state}

Provided $\mathbf{\big(C1\big)}$ holds, according to Corollary \ref{Corollary} the data $d_n(t)=x_n(t), \forall n\in{\cal{S}}$ in $\big(0, T^0\big)$ determines uniquely the network unknown initial conditions $X_0=X(0)$ and $\bar X_0=\dot X(0)$ in the problem (\ref{weq_0T0}). We write the solution 
$X$ of the problem (\ref{weq_0T0}) as follows:
\begin{eqnarray}
X(t)=X^{0}(t) + X^{F}(t), \;\; \forall t\in(0,T^0),
\label{splitXH}  
\end{eqnarray}
where $X^{0}=\big(x_1^0,\dots,x_N^0\big)^{\top}$ and $X^{F}=\big(x_1^F,\dots,x_N^F\big)^{\top}$ are the solutions of:
\begin{eqnarray}
\left\{
\begin{array}{lll}
\ddot X^{0}(t) + \eta \dot X^{0}(t)  - \Delta X^{0}(t) =0 \quad \mbox{in} \; (0,T^0),\\
X^{0}(0)=X_0\in\R^N  \quad \mbox{and} \quad \dot X^{0}(0)=\bar X_0\in\R^N,
\end{array}
\right.
\label{Healthy0_0T0}
\end{eqnarray}
and
\begin{eqnarray}
\left\{
\begin{array}{lll}
\ddot X^{F}(t) + \eta \dot X^{F}(t)  - \Delta X^{F}(t) =F^{sour}(t) \quad \mbox{in} \; (0,T^0),\\
X^{F}(0)=0  \quad \mbox{and} \quad \dot X^{F}(0)=0.
\end{array}
\right.
\label{HealthyF_0T0}
\end{eqnarray}
We expand the unknown initial conditions of the problem (\ref{Healthy0_0T0}) in the orthonormal basis of $\R^N$ defined by the normalized eigenvectors $v^{k_{\ell}}$ of the graph Laplacian matrix $\Delta$: 
\begin{eqnarray}
X_0=\displaystyle\sum_{k=1}^K\displaystyle\sum_{\ell=1}^{m_k} y^0_{k_{\ell}} v^{k_{\ell}}  \qquad \mbox{and} \qquad  \bar X_0=\displaystyle\sum_{k=1}^K\displaystyle\sum_{\ell=1}^{m_k} \bar y^0_{k_{\ell}} v^{k_{\ell}}.
\label{Expand_initialVectors_Delta}
\end{eqnarray}
Then, we establish the following result:
\begin{proposition}
Provided (\ref{Expand_initialVectors_Delta}) holds, the solution of the problem (\ref{Healthy0_0T0}) is given by
\begin{eqnarray}
X^{0}(t)=\displaystyle\sum_{k=1}^K\displaystyle\sum_{\ell=1}^{m_k} \Big({\cal{A}}_k(t) y^0_{k_{\ell}} + {\cal{B}}_k(t) \bar y^0_{k_{\ell}} \Big) v^{k_{\ell}}, \;\; \forall t\in(0,T^0),
\label{Health0_State_X0}
\end{eqnarray}
where the two functions ${\cal{A}}_k$ and ${\cal{B}}_k$ are defined as follows:
\begin{eqnarray}
{\cal{A}}_{k}(t)=
\left\{
\begin{array}{llll}
\displaystyle\frac{\bar r_k e^{r_kt} - r_k e^{\bar r_kt}}{\bar r_k-r_k}, \qquad\qquad\qquad\qquad\qquad\qquad\quad \mbox{if}\; D_k >0,\\
\big(1 + \frac{\eta}{2}t\big)e^{-\frac{\eta}{2}t}, \qquad\qquad\qquad\qquad\qquad\qquad\quad\;\; \mbox{if}\; D_k=0,\\
e^{-\frac{\eta}{2}t}\Big(\cos\big(\frac{\sqrt{-D_k}}{2}t\big) + \frac{\eta}{\sqrt{-D_k}}\sin\big(\frac{\sqrt{-D_k}}{2}t\big)\Big),  \qquad  \mbox{if}\; D_k<0,
\end{array}
\right.
\label{A}
\end{eqnarray}
\begin{eqnarray}
\hspace{-4cm}{\cal{B}}_{k}(t)=
\left\{
\begin{array}{llll}
\displaystyle\frac{e^{\bar r_kt} - e^{r_kt}}{\bar r_k-r_k}, \qquad\qquad\;\; \mbox{if}\; D_k >0,\\
t e^{-\frac{\eta}{2}t}, \qquad\qquad\qquad\;\; \mbox{if}\; D_k=0,\\
\frac{2 e^{-\frac{\eta}{2}t}}{\sqrt{-D_k}}\sin\big(\frac{\sqrt{-D_k}}{2}t\big),  \quad\;\;  \mbox{if}\; D_k<0,
\end{array}
\right.
\label{B}
\end{eqnarray}
with the discriminant $D_k=\eta^2 - 4\omega_k^2$ and the roots $r_{k}=\frac{-\eta + \sqrt{D_k}}{2}$ and $\bar r_{k}=-\frac{\eta + \sqrt{D_k}}{2}$.
\label{Solution_X0}
\end{proposition}  
{\bf{Proof.}} See the appendix.

According to (\ref{Health0_State_X0}), the identification of the unknown network initial conditions defining through (\ref{splitXH}) the state $X$ solution of the problem (\ref{weq_0T0}) can be achieved by determining the coefficients $y^0_{k_{\ell}}$ and $\bar y^0_{k_{\ell}}$ in (\ref{Expand_initialVectors_Delta}). To this end, we solve the minimization problem:
\begin{eqnarray}
\begin{array}{cccc}
  \displaystyle\min_{y^0_{k_{\ell}},\bar y^0_{k_{\ell}}}J\big(y^0_{k_{\ell}},\bar y^0_{k_{\ell}}\big)&:=&\frac{1}{2}\displaystyle\sum_{n\in{\cal{S}}}\big\|x_n(t) -d_n(t)\big\|^2_{L^2(0,T^0)},\qquad\qquad\qquad\qquad\qquad\qquad\qquad\qquad&\\
  &=&\frac{1}{2}\displaystyle\sum_{n\in{\cal{S}}}\int_0^{T^0}\Big(\displaystyle\sum_{k=1}^K\displaystyle\sum_{\ell=1}^{m_k} \Big({\cal{A}}_k(t) y^0_{k_{\ell}} + {\cal{B}}_k(t) \bar y^0_{k_{\ell}} \Big) v_n^{k_{\ell}} + x_n^{F}(t) - d_n(t)\Big)^2dt,&
\end{array}
\label{Min_J}
\end{eqnarray}
where from (\ref{splitXH}), $x_n(t)=x_n^0(t) + x_n^F(t), \forall n$.

We can write the gradient of the objective 
function $J$ in (\ref{Min_J}) in the following matrix form:
\begin{eqnarray}
\nabla J\big(Y\big)=M Y + b.
\label{nablaJ}
\end{eqnarray}
where $Y_{k+\ell-1}=y^0_{k_{\ell}},~~
Y_{N+k+\ell-1}=\bar y^0_{k_{\ell}}, ~~k=1,\dots K,~\ell=1,\dots m_k$. 
The expressions of the matrix $M$ 
and the right hand side $b$ are given in the Appendix with more details.

Writing $\nabla J=0$, we get the linear system $MY=-b$
whose solution gives the network initial condition as follows:
\begin{eqnarray}
X_0=\displaystyle\sum_{k=1}^K\displaystyle\sum_{\ell=1}^{m_k} Y_{k+\ell-1} v^{k_{\ell}} \qquad \mbox{and} \qquad  \bar X_0=\displaystyle\sum_{k=1}^K\displaystyle\sum_{\ell=1}^{m_k} Y_{N+k+\ell-1} v^{k_{\ell}}.
\label{Expand_initialVectors2}
\end{eqnarray}
Afterwards, from solving the forward problem (\ref{weq_0T0}) with the initial conditions obtained in (\ref{Expand_initialVectors2}) and $T$ instead of $T^0$, we determine the healthy network state  $X^H$ in $(0, T)$.

\subsection{Inverse problem for dominantly absorbent vertices}

We consider the network shown in Fig. \ref{g5}.
From the network shown in Fig. 2, we select the strategic and dominantly absorbent set of observation vertices ${\cal{S}} = \{1,4,5\}$. 
That leaves the set of non-observation vertices ${\cal{E}} = \{2,3\}$.
Then, we set $F^{sour}(t)=0$ and generate synthetic measurements 
from disturbances located in the vertex $m=2$ of intensity 
$$F^{dis}(t) = \sin {\pi(t-T_0) \over T},~~t> T_0, ~~0, ~{\rm otherwise} $$ 
where $T=100, ~T_0=10, ~\eta=0.1, ~\Delta t=0.1$. The initial 
condition of the system is $X_0= \bar X_0=0$.
\begin{figure} [H]
\centerline{
\epsfig{file=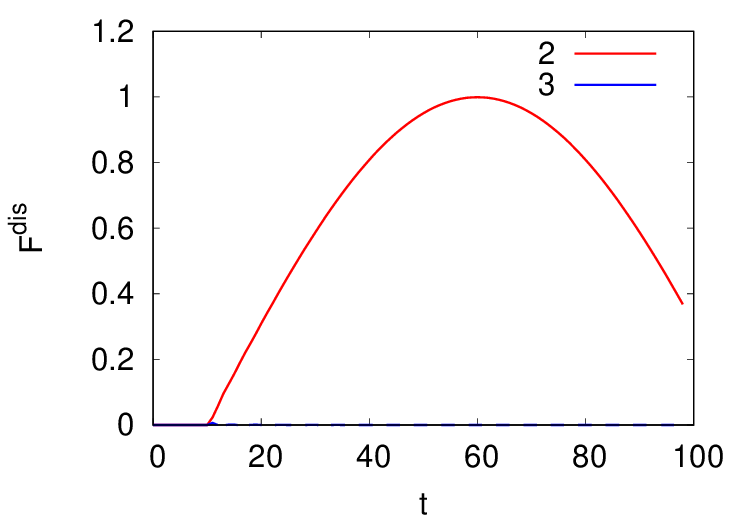,height=4 cm,width=6 cm,angle=0}
\epsfig{file=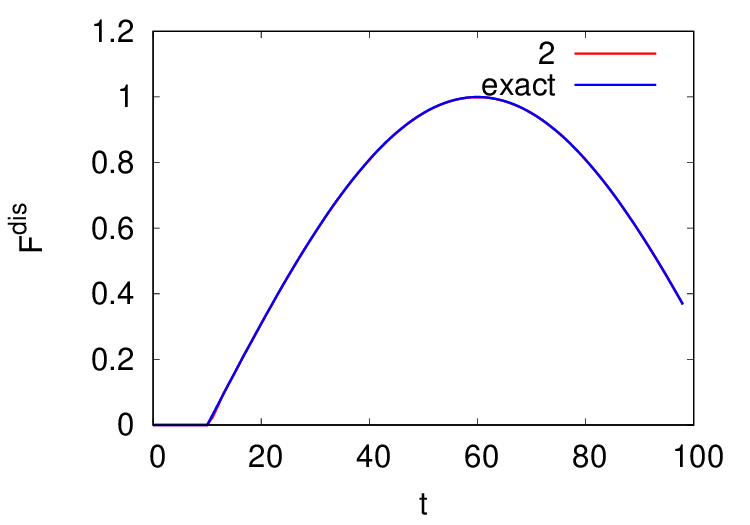,height=4 cm,width=6 cm,angle=0}
}
\caption{
Left panel: reconstructed disturbances in vertices 2 and 3,
Right panel: reconstructed and exact disturbances in vertex 2.
}
\label{sd2}
\end{figure}
Fig. \ref{sd2} shows the reconstructed disturbances in vertices 2 
and 3 (left panel) and the reconstructed and exact disturbances in 
vertices 2 and 3 (right panel). In the left panel, the reconstructed
disturbance at vertex 3 is zero as expected. 
The disturbance at vertex 2 is non zero and
is accurately reconstructed as seen in the right panel of Fig. \ref{sd2}.
Indeed, the reconstructed disturbances are non-null only in the vertex $m=2$ 
origin of these disturbances. Moreover, the reconstruction is
accurate because the unmeasured residuals in the vertices ${\cal{E}}$
are determined by solving the well posed linear system 
(\ref{Residual_LS_matrix}), 
where the full rank matrix is
$$\Delta\big({\cal{S}};{\cal{E}}\big) =
\begin{pmatrix}
0 & 0 \\
1 & 1 \\
1 & 0 
\end{pmatrix}.$$

\subsection{Inverse problem in absorbent observation vertices}

Now assume that the observation set is
${\cal{S}} = \{1,4\}$, it is not anymore dominantly absorbent
but only absorbent. The non observed vertices are
$ {\cal{E}} = \{2,3,5\}$. The reconstruction of the
residuals $x_m^R(t), ~m=2,3,5$ is done using the minimization
problem (\ref{Min_xi}) with a regularization factor $\alpha = 10^{-2}$.
\begin{figure} [H]
\centerline{
\epsfig{file=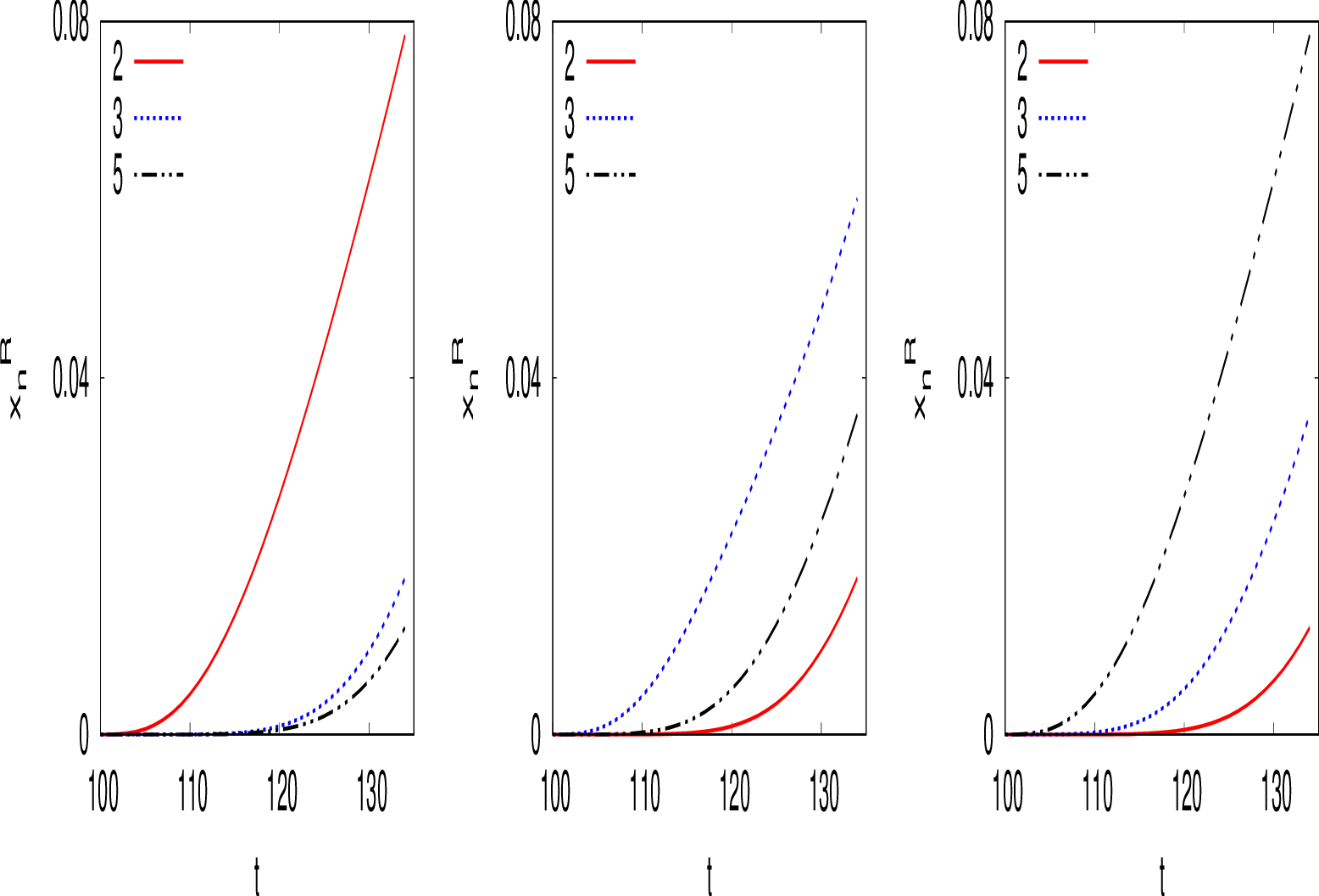,height=15 cm,width=6 cm,angle=90}
}
\caption{Reconstructed Residuals $x_n^R(t)$ for disturbances 
located in vertex: 2(left), 3(middle) and 5(right).}
\label{arr}
\end{figure}
Fig. \ref{arr} shows the time evolution of the residuals $x_m^R(t)$ for
a source term in vertices 2,3 or 5. 
These numerical results back our 
interpretations of the assertions (\ref{interpretation}). Indeed, 
these results show that the residual in the vertex origin of disturbances 
is the first to become significantly non-zero and is clearly larger
than the residuals in the other vertices wherein there are no 
disturbances. These elements define our procedure to localize 
disturbances using a non-dominantly absorbent observation set ${\cal{S}}$.
Notice that overtime, disturbances end up contaminating all the
vertices of the network so that residuals cannot be used anymore to
find the initial location and intensities of the detected disturbances.
Therefore, we reconstruct residuals and disturbances 
in the interval $(T_0, T_k)$, where $T_k$ is reasonably larger 
than the first instant of detection of disturbances.

\begin{figure} [H]
\centerline{
\epsfig{file=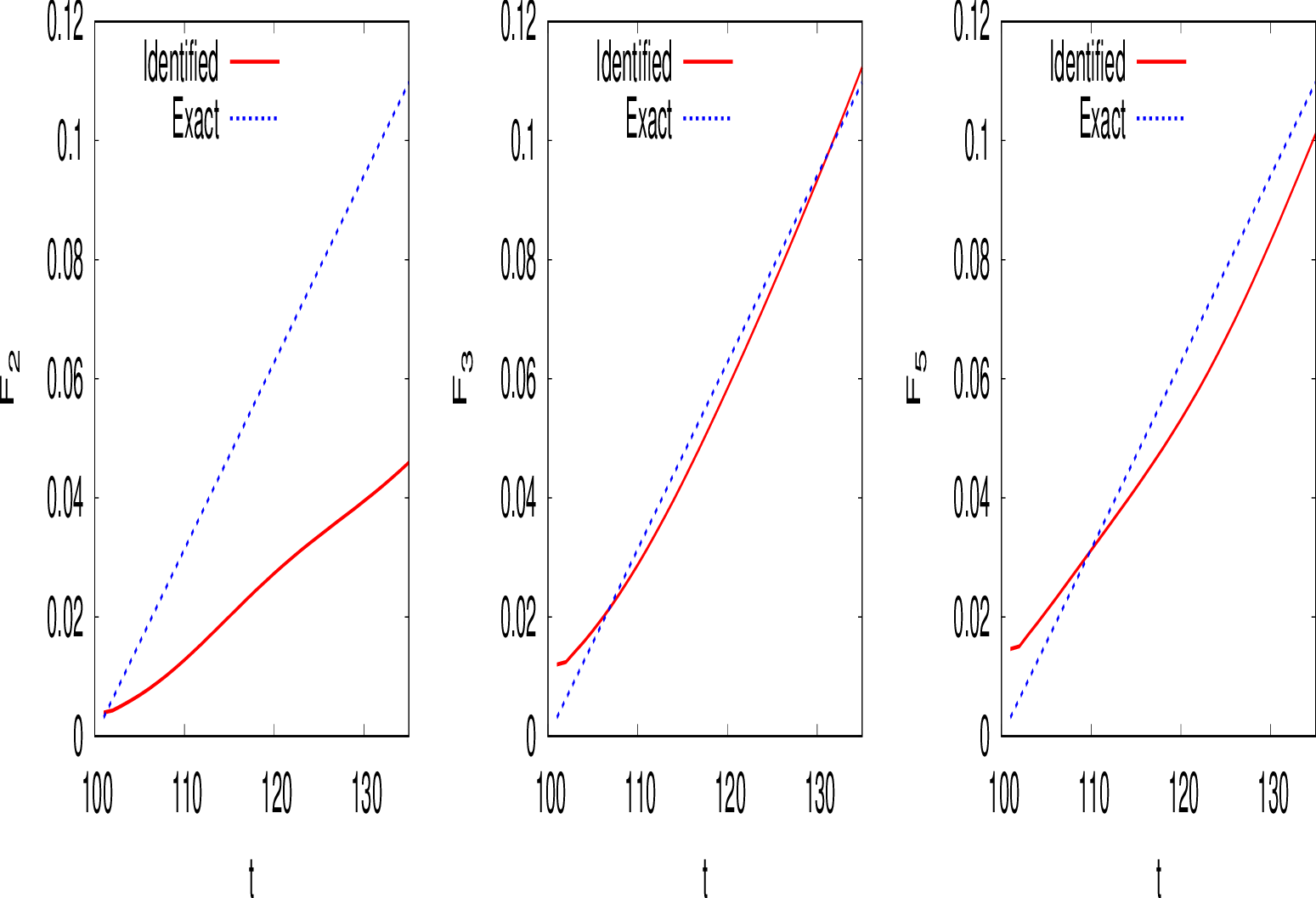,height=15 cm,width=6 cm,angle=90}
}
\caption{Reconstruction of disturbances $F_m^{dis}(t)$ located in vertex: 2(left), 3(middle) and 5(right).}
\label{adis}
\end{figure}

Fig. \ref{adis} shows the reconstruction of the disturbance $F_m^{dis}(t)$
once the source node has been found.
Note that $F_m^{dis}(t)$ is better approximated for a source in vertices
3 and 5 than in vertex 2. 
Two factors could amplify the error on the residual in a vertex $m$ 
of ${\cal{E}}$: the total number of network vertices directly 
connected to the vertex $m$ and how many of these are not observed. Since 
the residual $x_m^R$ is the major contributor in the computed data 
from (\ref{data_own_scheme}), this perspective might explain the 
errors on the reconstructed disturbances $F_m^{dis}$ presented in 
Fig. \ref{adis}. Indeed, the smallest error is on $F_3^{dis}$ since from 
Fig. \ref{g5} vertex $m=3$ is directly connected to 2 vertices 
with only one of them being not observed. As expected, the largest 
error is for $F_2^{dis}$, because vertex $m=2$ is directly connected 
to 3 vertices, 2 of them being not observed. The error on 
$F_5^{dis}$ is smaller than the one for $F_2^{dis}$ since only 1 
of the 3 directly connected vertices to $m=5$ is not observed.
On the other hand, it is larger than the error on $F_3^{dis}$ 
since $m=5$ is directly connected to more vertices than $m=3$.

%
\section{Conclusion}

In the present study, we addressed the inverse source problem of 
detection-localization-identification of unknown disturbances that 
might accidentally affect transmission networks after an initial 
healthy period of time.
This study is a follow up 
on our previous work where we broaden the spectrum of applications to 
cover more general networks defined by graph Laplacian matrices 
admitting eigenvalues of arbitrary multiplicities. We also 
allow a damping coefficient on the transmission process and 
consider that disturbances might affect multiple non-accessible 
vertices.

We generalized the notion of strategic set of observation vertices and 
proved that recording the network state in such a set enables 
to uniquely identify its initial state during the healthy period 
of time. From this we determine the global healthy network state.
This answers the question of the introduction on how to choose
the observation set to detect the presence of a disturbance.
A strategic observation set is enough to detect the
presence of a disturbance; we can use algorithm 1: detection of disturbances.

Then, we introduced the notion of dominantly absorbent set of 
network vertices. We proved that if the observation set is strategic 
and dominantly absorbent, then these observations determine in 
a unique manner the states of the network in the non-observation 
vertices since these solve a well posed linear system. This enables 
to efficiently detect and uniquely identify the unknown disturbances 
in quasi-real-time. This answers the third question of
the introduction: a strategic and dominantly absorbent observation
set allows quasi-real-time localization and identification of a 
disturbance. In practice, this can be done using Algorithm 2.

In practice, the number of observation vertices is limited and 
the observation set is not dominantly absorbent. Then
the linear system becomes ill-posed. To identify the
network states in the non-observation vertices, we need the observation
set to be at least absorbent. To obtain these states,
we introduce a least squares problem augmented by an 
appropriate regularization term.
The deviation 
of these states from their healthy values enables to localize the 
unknown disturbances and identify their intensities.
This answers question two of the introduction: we use in practice
Algorithm 3.


\newpage

\section{Appendix: }
\subsection{An algorithm to find an absorbent set of a graph}

A possible algorithm to find an absorbing set is to 
find a covering tree of the graph using for example Kurskal's algorithm,
see \cite{precis}. Then one decimates the vertices
of largest degree until the remaining set is absorbent.

\hspace{0cm}{\bf{\uwave{Algorithm}: Absorbing set}} \\
{\bf{Choose:}} ${\cal{S}}$ empty set

\hspace{0.2cm}{\bf{\underline{Begin}}}
\begin{itemize}

\item [$\bullet$]
Find a covering tree of the graph using Kruskal's algorithm. 
Save the vertices in the set ${\cal{T}}$.

\item [$\bullet$] Find maximal degree $d_{max}$ of 
the vertices in ${\cal{T}}$.

\item [$\bullet$]
~~For $d = d_{max}$ to 2 \\
Select all vertices in ${\cal{T}}$ of degree $ \ge d$  \\
Include them in ${\cal{S}}$ \\
If ${\cal{S}}$ is not absorbent then Break \\
~~end For
\end{itemize}
{\bf{\underline{End}.}}

\subsection{Proof of Theorem \ref{Th_InitialConditions}}                             

Let $X_0\in\R^N$, $\bar X_0\in\R^N$ and $X^{R}=\big(x^{R}_1,\dots,x^{R}_N\big)^{\top}$ be the solution of the problem:
\begin{eqnarray}
\left\{
\begin{array}{lll}
\ddot X^{R}(t) + \eta \dot X^{R}(t)  - \Delta X^{R}(t) =0 \quad \mbox{in} \; \big(0,T^0\big),\\
X^{R}(0)=X_0  \quad \mbox{and} \quad \dot X^{R}(0)=\bar X_0.
\end{array}
\right.
\label{Rweq_0T0_Proof}
\end{eqnarray}
From expanding the solution $X^R$ of the problem (\ref{Rweq_0T0_Proof}) in the orthonormal basis of $\R^N$ defined by the normalized eigenvectors of the graph Laplacian matrix $\Delta$, we write
\begin{eqnarray}
X^{R}(t)=\displaystyle\sum_{k=1}^K\displaystyle\sum_{\ell=1}^{m_k} y_{k_{\ell}}(t) v^{k_{\ell}}, \;\; \forall t, \qquad \mbox{where}  \;\; y_{k_{\ell}}(t)=\langle X^{R}(t),v^{k_{\ell}}\rangle.
\label{Global_Health_State}
\end{eqnarray}
Here, $\langle,\rangle$ denotes the inner product in $\R^N$. Moreover, according to (\ref{eigenvectors}), the functions $y_{k_{\ell}}$ in (\ref{Global_Health_State}) solve the following second order differential equation:
\begin{eqnarray}
\left\{
\begin{array}{llll}  
\ddot y_{k_{\ell}}(t) + \eta \dot y_{k_{\ell}}(t) + \omega^2_{k}y_{k_{\ell}}(t)=0 \quad \mbox{in} \; \big(0,T^0\big),\\
y_{k_{\ell}}(0)=\langle X_0,v^{k_{\ell}}\rangle \quad \mbox{and} \quad \dot y_{k_{\ell}}(0)=\langle \bar X_0,v^{k_{\ell}}\rangle.
\end{array}
\right.
\label{Diff_Eq}
\end{eqnarray}
The discriminant of the characteristic equation associated to (\ref{Diff_Eq}) is $D_k=\eta^2 - 4\omega_k^2$. Depending on the sign of $D_k$, the general form of the real-valued solutions of (\ref{Diff_Eq}) reads
\begin{eqnarray}
y_{k_{\ell}}(t)=e^{-\frac{\eta}{2}t}
\left\{
\begin{array}{llll}
\lambda_{k_{\ell}} e^{\frac{\sqrt{D_k}}{2}t}  + \mu_{k_{\ell}} e^{-\frac{\sqrt{D_k}}{2}t}, \qquad\qquad\qquad\quad \mbox{if}\; D_k >0,\\
\lambda_{k_{\ell}} + \mu_{k_{\ell}} t, \qquad\qquad\qquad\qquad\qquad\qquad \mbox{if}\; D_k=0,\\
\lambda_{k_{\ell}}\cos\big(\frac{\sqrt{-D_k}}{2}t\big) + \mu_{k_{\ell}}\sin\big(\frac{\sqrt{-D_k}}{2}t\big),  \qquad  \mbox{if}\; D_k<0,
\end{array}
\right.
\label{y_global}
\end{eqnarray}
where the real coefficients $\lambda_{k_{\ell}}$ and $\mu_{k_{\ell}}$ to be determined from the initial conditions $y_{k_{\ell}}(0)$ and $\dot y_{k_{\ell}}(0)$. Assuming ${\cal{S}}$ to be a {\it{strategic}} set of vertices, suppose that
\begin{eqnarray}
 x_n^{R}(t)=0, \;\; \forall t\in \big(0,T^0\big),\; \forall n\in {\cal{S}}.
\label{Hypothesis}
\end{eqnarray}
Then, in view of (\ref{Global_Health_State})-(\ref{Hypothesis}), it follows that
\begin{eqnarray}
x_n^{R}(t)=\displaystyle\sum_{k=1}^K\displaystyle\sum_{\ell=1}^{m_k} y_{k_{\ell}}(t) v_n^{k_{\ell}}=0, \;\; \forall t\in \big(0,T^0\big), \; \forall n \in {\cal{S}}.
\label{assumption}
\end{eqnarray}
Besides for arbitrary values $\eta>0$ and $\omega_k^2\ge 0$, we notice that according to (\ref{eigenValues}) and since $\omega_1=0$, the set of discriminants fulfills
\begin{eqnarray}
\eta^2-4\omega_K^2=D_K<D_{K-1}<\dots <D_2<D_1=\eta^2. 
\label{Discriminant_order}
\end{eqnarray}
To establish the proof for the most general case of the three possible forms defining the functions $y_{k_{\ell}}$ in (\ref{y_global}) are involved in (\ref{assumption}), we assume that
\begin{eqnarray}
\exists \bar k\in\{2,\dots,K-1\} \quad \mbox{such that} \quad D_{\bar k}=0. 
\label{Discriminant}
\end{eqnarray}
According to (\ref{y_global})-(\ref{Discriminant}), it follows from (\ref{assumption}) by analytic continuation that: $\forall n \in {\cal{S}}$,
\begin{eqnarray}
\begin{array}{llll}
  \displaystyle\sum_{k=1}^{\bar k-1}\displaystyle\sum_{\ell=1}^{m_k}\Big(\lambda_{k_{\ell}} e^{\frac{\sqrt{D_k}}{2}t}  + \mu_{k_{\ell}} e^{-\frac{\sqrt{D_k}}{2}t}\Big)v_n^{k_{\ell}} + \displaystyle\sum_{\ell=1}^{m_{\bar k}}\Big(\lambda_{\bar k_{\ell}} + \mu_{\bar k_{\ell}}t\Big)v_n^{\bar k_{\ell}}\\
+ \displaystyle\sum_{k=\bar k + 1}^{K}\displaystyle\sum_{\ell=1}^{m_k}\Big(\lambda_{k_{\ell}}\cos\big(\frac{\sqrt{-D_k}}{2}t\big) + \mu_{k_{\ell}}\sin\big(\frac{\sqrt{-D_k}}{2}t\big)\Big)v_n^{k_{\ell}}=0, \;\; \forall t\in\R.
\end{array}
\label{LocalState}
\end{eqnarray}
Notice that: {\bf{1.}} If $\bar k$ doesn't exist and $D_K>0$, then the Left Hand Side (LHS) of (\ref{LocalState}) is defined only by its first term with $K$ instead of $\bar k-1$. {\bf{2.}} If $\bar k$ doesn't exist and $D_K<0$, then the LHS of (\ref{LocalState}) is defined only by its first and third terms with $\hat k$ instead of $\bar k-1$ and $\hat k + 1$ instead of $\bar k + 1$, where $\hat k$ is such that $D_{\hat k}>0$ whereas $D_{\hat k + 1}<0$. {\bf{3.}} If $\bar k=K$, then $D_K=0$ and the LHS of (\ref{LocalState}) is defined only by its two first terms.

For $p=1,\dots,K$, let $A^{(p)}$ be the $N_{\cal{S}}\times m_p$ matrix whose columns are the $m_p$ eigenvectors of $\Delta$ associated to $-\omega_p^2$ and rows are the $N_{\cal{S}}$ components of these vectors corresponding to the vertex labels defining ${\cal{S}}$. Since ${\cal{S}}$ is strategic, it follows from Definition \ref{strategic_set} that the matrix $A^{(p)}$ is of full rank. Let $\Lambda^{(p)}=\big(\lambda_{p_{1}}, \dots, \lambda_{p_{m_p}}\big)^{\top}$ and $\zeta^{(p)}=\big(\mu_{p_{1}}, \dots, \mu_{p_{m_p}}\big)^{\top}$.

\hspace{-0.3cm} $\blacktriangleright$ {\bf{Proving that $\lambda_{k_{\ell}}=\mu_{k_{\ell}}=0$, for $k=1,\dots,\bar k-1$ and $\ell=1,\dots,m_k$:}}

\hspace{-0.1cm}$\bullet$ {\it{Using $D_{p=1}$:}} Since $D_1>D_{k=2,\dots,K}$, from multiplying (\ref{LocalState}) by $e^{-{\frac{\sqrt{D_1}}{2}t}}$ and setting the limit when $t$ tends to $+\infty$, then by $e^{{\frac{\sqrt{D_1}}{2}t}}$ and setting the limit when $t$ goes to $-\infty$, we get
\begin{eqnarray}
\forall n \in {\cal{S}}, \qquad \displaystyle\sum_{\ell=1}^{m_1} \lambda_{1_{\ell}} v_n^{1_{\ell}}=0 \qquad \mbox{and} \qquad \displaystyle\sum_{\ell=1}^{m_1} \mu_{1_{\ell}} v_n^{1_{\ell}}=0.
\label{identification_k1}
\end{eqnarray}
Moreover, (\ref{identification_k1}) can be rewritten under the following matrix form:
\begin{eqnarray}
A^{(1)}\Lambda^{(1)}=0 \quad \mbox{and}  \quad A^{(1)}\zeta^{(1)}=0  \quad \implies \quad \lambda_{1_{\ell=1,\dots,m_1}}=\mu_{1_{\ell=1,\dots,m_1}}=0.
\label{Final-identification_k1}
\end{eqnarray}
The implication in (\ref{Final-identification_k1}) is obtained since the matrix $A^{(1)}$ is of full rank.

$\bullet$ {\it{Using $D_{p=2,\dots,\bar k-1}$:}} By taking into account of (\ref{Final-identification_k1}) in (\ref{LocalState}) and iterating the same process for $p=2,\dots,\bar k-1$: We multiply (\ref{LocalState}) by $e^{-{\frac{\sqrt{D_p}}{2}t}}$ and set the limit when $t$ tends to $+\infty$, then by $e^{{\frac{\sqrt{D_p}}{2}t}}$ and set the limit when $t$ goes to $-\infty$, we obtain: For $p=2,\dots,\bar k-1$,
\begin{eqnarray}
A^{(p)}\Lambda^{(p)}=0 \quad \mbox{and}  \quad A^{(p)}\zeta^{(p)}=0  \quad \implies \quad \lambda_{p_{\ell=1,\dots,m_p}}=\mu_{p_{\ell=1,\dots,m_p}}=0.
\label{Final-identification_k1_bark}
\end{eqnarray}
\hspace{-0.3cm} $\blacktriangleright$ {\bf{Proving that $\lambda_{\bar k_{\ell}}=\mu_{\bar k_{\ell}}=0$, for $\ell=1,\dots,m_{\bar k}$:}} In view of (\ref{Final-identification_k1})-(\ref{Final-identification_k1_bark}), the first double sum in (\ref{LocalState}) vanishes. We compute the integral over $(0,{\cal{T}})$ of the updated equation (\ref{LocalState}), where ${\cal{T}}>0$. From dividing the obtained result, firstly by ${\cal{T}}^2$ and setting the limit when ${\cal{T}}$ tends to $+\infty$, then by ${\cal{T}}$ and resetting again the limit when ${\cal{T}}$ goes to $+\infty$, we find
\begin{eqnarray}
\forall n \in {\cal{S}}, \quad \displaystyle\sum_{\ell=1}^{m_{\bar k}} \mu_{\bar k_{\ell}} v_n^{\bar k_{\ell}}=0 \quad \mbox{and} \quad \displaystyle\sum_{\ell=1}^{m_{\bar k}} \lambda_{\bar k_{\ell}} v_n^{\bar k_{\ell}}=0 \quad \implies \quad \lambda_{\bar k_{\ell=1,\dots,m_{\bar k}}}=\mu_{\bar k_{\ell=1,\dots,m_{\bar k}}}=0.
\label{Final-identification_bark}
\end{eqnarray}
The implication in (\ref{Final-identification_bark}) is obtained since the matrix $A^{(\bar k)}$ is of full rank.

\hspace{-0.3cm} $\blacktriangleright$ {\bf{Proving that $\lambda_{k_{\ell}}=\mu_{k_{\ell}}=0$, for $k=\bar k+1,\dots,K$ and $\ell=1,\dots,m_k$:}} According to (\ref{Final-identification_k1_bark})-(\ref{Final-identification_bark}), it follows that (\ref{LocalState}) is now reduced to:
\begin{eqnarray}
\forall n \in {\cal{S}}, \quad \displaystyle\sum_{k=\bar k + 1}^{K}\displaystyle\sum_{\ell=1}^{m_k}\Big(\lambda_{k_{\ell}}\cos\big(\frac{\sqrt{-D_k}}{2}t\big) + \mu_{k_{\ell}}\sin\big(\frac{\sqrt{-D_k}}{2}t\big)\Big)v_n^{k_{\ell}}=0, \;\; \forall t\in\R.
\label{LocalState_Reduced}
\end{eqnarray}
For $p=\bar k+1,\dots,K$: We multiply (\ref{LocalState_Reduced}) by $\cos\big(\frac{\sqrt{-D_p}}{2}t\big)$ and compute its integral over $(0,{\cal{T}})$. From dividing the result by ${\cal{T}}$ and setting the limit when ${\cal{T}}$ tends to $+\infty$, we get  
\begin{eqnarray}
\forall n \in {\cal{S}}, \quad \displaystyle\sum_{\ell=1}^{m_{p}} \lambda_{p_{\ell}} v_n^{p_{\ell}}=0 \quad \implies \quad \lambda_{p_{\ell=1,\dots,m_p}}=0. 
\label{identification_cos_p}
\end{eqnarray}
Afterwards, by taking into account the result (\ref{identification_cos_p}) in (\ref{LocalState_Reduced}) and then, computing the derivative of the updated equation, it follows that 
\begin{eqnarray}
\forall n \in {\cal{S}}, \quad \displaystyle\sum_{k=\bar k + 1}^{K}\displaystyle\sum_{\ell=1}^{m_k}\frac{\sqrt{-D_k}}{2}\mu_{k_{\ell}}\cos\big(\frac{\sqrt{-D_k}}{2}t\big) v_n^{k_{\ell}}=0, \quad \forall t\in\R.
\label{LocalState_Reduced2}
\end{eqnarray}
Then, reapplying the same techniques used in (\ref{LocalState_Reduced})-(\ref{identification_cos_p}), we obtain: For $p=\bar k+1,\dots,K$, 
\begin{eqnarray}
\forall n \in {\cal{S}}, \quad \displaystyle\sum_{\ell=1}^{m_{p}} \mu_{p_{\ell}} v_n^{p_{\ell}}=0 \quad \implies \quad \mu_{p_{\ell=1,\dots,m_p}}=0. 
\label{identification_sin_p}
\end{eqnarray}
The implications in (\ref{identification_cos_p}) and in (\ref{identification_sin_p}) are obtained since the matrix $A^{(p)}$ is of full rank.

In view of (\ref{Final-identification_k1})-(\ref{Final-identification_k1_bark}), (\ref{Final-identification_bark}), (\ref{identification_cos_p}) and (\ref{identification_sin_p}), it follows that $\lambda_{k_{\ell=1,\dots,m_k}}=\mu_{k_{\ell=1,\dots,m_k}}=0$, for $k=1,\dots,K$. Therefore, (\ref{Global_Health_State})-(\ref{y_global}) implies that the solution of the problem (\ref{Rweq_0T0_Proof}) is $X^R(t)=0, \forall t$ which leads to $X^R(0)=\dot X^R(0)=0$. \hspace{6cm} $\blacksquare$

\subsection{Proof of Propostion \ref{Own_numerical_scheme}}

Let $a<b$ be two real numbers, $t_i=(a+b)/2$ and $\varphi_i$ be the function introduced in (\ref{Approx_phi}) with $\big(a, b\big)$ instead of $\big(t_{i-1}, t_{i+1}\big)$. Let $\alpha, \beta, \gamma$ be three real numbers and $g(t)=\alpha t^2 + \beta t + \gamma, \forall t$.
\begin{eqnarray*}
\begin{array}{ccccc}
\blacktriangleright  \displaystyle\int_a^{t_i}g(t)\varphi_i(t)dt&=&-\displaystyle\frac{1}{2}\int_a^{t_i}\Big(\alpha t^3 + \big[\beta - a\alpha\big]t^2 + \big[\gamma-\beta a\big]t - \gamma a\Big)dt, \qquad \qquad \qquad\qquad \qquad \qquad&\\
  &=&-\displaystyle\frac{1}{2}\Big(\frac{\alpha}{4}\big(t_i^4 - a^4\big) + \frac{1}{3}\big[\beta - a\alpha\big]\big(t_i^3 - a^3\big) + \frac{1}{2}\big[\gamma-\beta a\big]\big(t_i^2 - a^2\big) - \gamma a \big(t_i - a\big)\Big).&
\end{array}
\label{integ_a_ti}
\end{eqnarray*}
\begin{eqnarray*}
\begin{array}{ccccc}
\blacktriangleright  \displaystyle\int_{t_i}^bg(t)\varphi_i(t)dt&=&\displaystyle\frac{1}{2}\int_{t_i}^b\Big(\alpha t^3 + \big[\beta - b\alpha\big]t^2 + \big[\gamma-\beta b\big]t - \gamma b\Big)dt, \qquad \qquad \qquad\qquad \qquad \qquad&\\
  &=&\displaystyle\frac{1}{2}\Big(\frac{\alpha}{4}\big(b^4 - t_i^4\big) + \frac{1}{3}\big[\beta - b\alpha\big]\big(b^3 - t_i^3\big) + \frac{1}{2}\big[\gamma-\beta b\big]\big(b^2 - t_i^2\big) - \gamma b \big(b- t_i\big)\Big).&
\end{array}
\label{integ_ti_b}
\end{eqnarray*}
Since $t_i-a=b-t_i$, from employing Chasles rule it follows that
\begin{eqnarray}
\displaystyle\int_a^bg(t)\varphi_i(t)dt=-\displaystyle\frac{b-t_i}{2}\Big(\frac{\alpha}{12}\big(b^3 - a^3\big) + \big[\frac{\alpha}{12}t_i + \frac{\beta}{6}\big]\big(b^2 - a^2\big) + \big[\frac{\alpha}{12}t_i^2 + \frac{\beta}{6}t_i + \frac{\gamma}{2}\big]\big(b-a\big)\Big).   
\label{Chasles}
\end{eqnarray}
Besides, $t_i=(a + b)/2$ implies that $ab=2t_i^2 -(b^2 + a^2)/2$. Therefore, we get
\begin{eqnarray}
b^3-a^3=\big(b-a\big)\big(b^2 + ab + a^2\big)=\frac{b-a}{2}\big(b^2 + 4t_i^2 + a^2\big) \quad \mbox{and} \quad b^2 - a^2=2\big(b-a\big)t_i.
\label{ab}
\end{eqnarray}
Thus, from using (\ref{ab}) in (\ref{Chasles}) and since $2\beta t_i=\beta(a+b)$, we obtain
\begin{eqnarray}
\begin{array}{ccccc}
  \displaystyle\int_a^bg(t)\varphi_i(t)dt&=&-\displaystyle\frac{b-t_i}{2}\displaystyle\frac{b-a}{2}\Big(\frac{10}{12}\alpha t_i^2 + \beta t_i +  \frac{\alpha}{12}\big(b^2 + a^2\big) + \gamma\Big), \qquad\qquad\qquad&\\
  &=&-\displaystyle\frac{b-t_i}{2}\displaystyle\frac{b-a}{2}\Big(\frac{10}{12}\alpha t_i^2 + \frac{10}{12}\beta t_i + \frac{\beta}{12} (a + b) +  \frac{\alpha}{12}\big(b^2 + a^2\big) + \gamma\Big),&\\
&=&-\displaystyle\frac{b-t_i}{2}\displaystyle\frac{b-a}{2}\Big(\frac{10}{12}g(t_i) + \frac{1}{12}g(b) + \frac{1}{12}g(a)\Big).\qquad\qquad\qquad\qquad& \\  
\end{array}
\label{Chasles_ab}
\end{eqnarray}
Since $t_i-a=b-t_i$, from (\ref{Approx_phi}) it follows that $\varphi_i(t_i)=-(b-t_i)/2$. Hence, the last equality in (\ref{Chasles_ab}) yields the result annouced in (\ref{Own_integration}). \hspace{6cm} $\blacksquare$

\subsection{Proof of Proposition \ref{Solution_X0}}

We expand the solution $X^0$ of the problem (\ref{Healthy0_0T0}) in the orthonormal basis defined by the family $\{v^{k_{\ell}}, \; k=1,\dots,K \;\mbox{and}\; \ell=1,\dots,m_k\}$ of the $N$ normalized eigenvectors associated to the graph Laplacian matrix $\Delta$:
\begin{eqnarray}
X^{0}(t)=\displaystyle\sum_{k=1}^K\displaystyle\sum_{\ell=1}^{m_k} y_{k_{\ell}}(t) v^{k_{\ell}}, \;\; \forall t\in(0,T^0),
\label{Health0_State}
\end{eqnarray}
where $y_{k_{\ell}}(t)=\langle X^{0}(t),v_{k_{\ell}}\rangle$ with $\langle,\rangle$ designates the inner product in $\R^N$. From the problem (\ref{Healthy0_0T0}), it follows in view of (\ref{eigenvectors}) and (\ref{Expand_initialVectors_Delta}) that the functions $y_{k_{\ell}}$ solve:
\begin{eqnarray}
\left\{
\begin{array}{llll}  
\ddot y_{k_{\ell}}(t) + \eta \dot y_{k_{\ell}}(t) + \omega^2_{k}y_{k_{\ell}}(t)=0, \;\;\forall t\in(0,T^0),\\
y_{k_{\ell}}(0)=y^0_{k_{\ell}} \quad \mbox{and} \quad \dot y_{k_{\ell}}(0)=\bar y^0_{k_{\ell}}.
\end{array}
\right.
\label{Diff_y}
\end{eqnarray}
The discriminant of the characteristic equation associated to (\ref{Diff_y}) is $D_k=\eta^2 - 4\omega_k^2$. Therefore, we distinguish the following three cases defining the solutions of (\ref{Diff_y}):
\begin{eqnarray}
y_{k_{\ell}}(t)=e^{-\frac{\eta}{2}t}
\left\{
\begin{array}{llll}
\lambda_{k_{\ell}} e^{\frac{\sqrt{D_k}}{2}t}  + \mu_{k_{\ell}} e^{-\frac{\sqrt{D_k}}{2}t}, \qquad\qquad\qquad\quad \mbox{if}\; D_k >0,\\
\lambda_{k_{\ell}} + \mu_{k_{\ell}} t, \qquad\qquad\qquad\qquad\qquad\qquad \mbox{if}\; D_k=0,\\
\lambda_{k_{\ell}}\cos\big(\frac{\sqrt{-D_k}}{2}t\big) + \mu_{k_{\ell}}\sin\big(\frac{\sqrt{-D_k}}{2}t\big),  \qquad  \mbox{if}\; D_k<0,
\end{array}
\right.
\label{y_express}
\end{eqnarray}
where the real coefficients $\lambda_{k_{\ell}}$ and $\mu_{k_{\ell}}$ to be determined from the initial conditions $y_{k_{\ell}}(0)$ and $\dot y_{k_{\ell}}(0)$. Then, depending on the sign of the discriminant $D_k$, we obtain:\\

\hspace{-0.3cm} $\blacktriangleright$ {\bf{\underline{Case 1}.}} If $D_k>0$, we denote by: $r_{k}=\frac{-\eta + \sqrt{D_k}}{2}$ and $\bar r_{k}=-\frac{\eta + \sqrt{D_k}}{2}$. From (\ref{y_express}), we get
\begin{eqnarray}
\left\{    
\begin{array}{lll}
\lambda_{k_{\ell}} + \mu_{k_{\ell}}=y_{k_{\ell}}(0),\\
\lambda_{k_{\ell}} r_{k} +  \mu_{k_{\ell}} \bar r_{k}=\dot y_{k_{\ell}}(0).
\end{array}
\right.
\quad \implies \quad
\begin{array}{lll}
\lambda_{k_{\ell}}=\frac{\bar r_k}{\bar r_k-r_k} y_{k_{\ell}}(0) -  \frac{1}{\bar r_k -r_k}\dot y_{k_{\ell}}(0),\\
\mu_{k_{\ell}}=-\frac{r_k}{\bar r_k - r_k} y_{k_{\ell}}(0) + \frac{1}{\bar r_k -r_k}\dot y_{k_{\ell}}(0).
\end{array}
\label{Dk_pos}
\end{eqnarray}
Hence, in this case, the solution $y_{k_{\ell}}$ of (\ref{Diff_y}) is then given by: 
\begin{eqnarray}
y_{k_{\ell}}(t)=\displaystyle\frac{\bar r_k e^{r_kt} - r_k e^{\bar r_kt}}{\bar r_k-r_k}y_{k_{\ell}}(0) + \frac{e^{\bar r_kt} - e^{r_kt}}{\bar r_k-r_k}\dot y_{k_{\ell}}(0). 
\label{y_Dk_Pos}
\end{eqnarray}
\hspace{-0.3cm} $\blacktriangleright$ {\bf{\underline{Case 2}.}} If $D_k=0$, then from (\ref{y_express}), it follows that the solution $y_{k_{\ell}}$ of (\ref{Diff_y}) is defined by:
\begin{eqnarray}
y_{k_{\ell}}(t)=\big(1 + \frac{\eta}{2}t\big)e^{-\frac{\eta}{2}t}y_{k_{\ell}}(0) + t e^{-\frac{\eta}{2}t}\dot y_{k_{\ell}}(0). 
\label{y_Dk_null}
\end{eqnarray}
\hspace{-0.3cm} $\blacktriangleright$ {\bf{\underline{Case 3}.}} If $D_k<0$, then (\ref{y_express}) leads to:
\begin{eqnarray}
y_{k_{\ell}}(t)=e^{-\frac{\eta}{2}t}\Big(\cos\big(\frac{\sqrt{-D_k}}{2}t\big) + \frac{\eta}{\sqrt{-D_k}}\sin\big(\frac{\sqrt{-D_k}}{2}t\big)\Big)y_{k_{\ell}}(0) + \frac{2 e^{-\frac{\eta}{2}t}}{\sqrt{-D_k}}\sin\big(\frac{\sqrt{-D_k}}{2}t\big)\dot y_{k_{\ell}}(0). 
\label{y_Dk_neg}
\end{eqnarray}
Setting $y_{k_{\ell}}(0)=y^0_{k_{\ell}}$ and $\dot y_{k_{\ell}}(0)=\bar y^0_{k_{\ell}}$ in (\ref{y_Dk_Pos})-(\ref{y_Dk_neg}) leads to the result in Proposition \ref{Solution_X0}. \hspace{0cm} $\blacksquare$

\subsection{Details on the gradient $\nabla J$ (\ref{nablaJ}) }

We verify that the first order partial derivatives of the objective function $J$ in (\ref{Min_J}) are given by: Let ${\cal{V}}_{k_{\ell}}^{p_q}=\displaystyle\sum_{n\in{\cal{S}}}v_n^{p_q}v_n^{k_{\ell}}$,
\begin{eqnarray}
\begin{array}{cccc}
\bullet \quad \displaystyle\frac{\partial J\big(y^0_{k_{\ell}},\bar y^0_{k_{\ell}}\big)}{\partial y^0_{p_q}}&=&\displaystyle\sum_{k=1}^K\displaystyle\sum_{\ell=1}^{m_k}\left({\cal{V}}_{k_{\ell}}^{p_q}\int_0^{T^0}{\cal{A}}_{p}(t){\cal{A}}_k(t)dt\; y^0_{k_{\ell}} + {\cal{V}}_{k_{\ell}}^{p_q}\int_0^{T^0}{\cal{A}}_{p}(t){\cal{B}}_k(t)dt\; \bar y^0_{k_{\ell}}\right)&\\
&+&\displaystyle\sum_{n\in{\cal{S}}}v^n_{p_q}\int_0^{T^0}{\cal{A}}_{p}(t)\Big( x_n^{F}(t) - d_n(t)\Big)dt. \qquad \qquad\qquad \qquad\qquad \qquad&
\end{array}
\label{Form2_GradientJ_y0pq}
\end{eqnarray}
\begin{eqnarray}
\begin{array}{cccc}
\bullet \quad \displaystyle\frac{\partial J\big(y^0_{k_{\ell}},\bar y^0_{k_{\ell}}\big)}{\partial \bar y^0_{p_q}}&=&\displaystyle\sum_{k=1}^K\displaystyle\sum_{\ell=1}^{m_k}\left({\cal{V}}_{k_{\ell}}^{p_q}\int_0^{T^0}{\cal{B}}_{p}(t){\cal{A}}_k(t)dt\; y^0_{k_{\ell}} + {\cal{V}}_{k_{\ell}}^{p_q}\int_0^{T^0}{\cal{B}}_{p}(t){\cal{B}}_k(t)dt\; \bar y^0_{k_{\ell}}\right)&\\
&+&\displaystyle\sum_{n\in{\cal{S}}}v_n^{p_q}\int_0^{T^0}{\cal{B}}_{p}(t)\Big( x_n^{F}(t) - d_n(t)\Big)dt. \qquad \qquad\qquad \qquad\qquad \qquad&
\end{array}
\label{Form2_GradientJ_dot_y0pq}
\end{eqnarray}
To express (\ref{Form2_GradientJ_y0pq})-(\ref{Form2_GradientJ_dot_y0pq}) under a matrix form, we introduce the $2N\times 2N$ matrix $M$ and the two vectors $b,Y$ of $\R^{2N}$ defined by: For $\left\{\begin{array}{lll}p=1,\dots,K\\ q=1,\dots,m_p\end{array}\right.\;$ and $\;\; \left\{\begin{array}{lll}k=1,\dots,K\\ \ell=1,\dots,m_k\end{array}\right.$,
\begin{eqnarray}
\left\{
\begin{array}{lllll}
M_{p+q-1,k+\ell-1}={\cal{V}}_{k_{\ell}}^{p_q}\displaystyle\int_0^{T^0}{\cal{A}}_{p}(t){\cal{A}}_k(t)dt, \qquad\quad M_{p+q-1,N+k+\ell-1}={\cal{V}}_{k_{\ell}}^{p_q}\displaystyle\int_0^{T^0}{\cal{A}}_{p}(t){\cal{B}}_k(t)dt,\\
M_{N+p+q-1,k+\ell-1}={\cal{V}}_{k_{\ell}}^{p_q}\displaystyle\int_0^{T^0}{\cal{B}}_{p}(t){\cal{A}}_k(t)dt, \qquad M_{N+p+q-1,N+k+\ell-1}={\cal{V}}_{k_{\ell}}^{p_q}\displaystyle\int_0^{T^0}{\cal{B}}_{p}(t){\cal{B}}_k(t)dt.
\end{array}
\right.
\label{M}
\end{eqnarray}
\begin{eqnarray}
\left\{
\begin{array}{llll} 
b_{p+q-1}=\displaystyle\sum_{n\in{\cal{S}}}v_n^{p_q}\int_0^{T^0}{\cal{A}}_{p}(t)\Big( x_n^{F}(t) - d_n(t)\Big)dt,\\
b_{N+p+q-1}=\displaystyle\sum_{n\in{\cal{S}}}v_n^{p_q}\int_0^{T^0}{\cal{B}}_{p}(t)\Big( x_n^{F}(t) - d_n(t)\Big)dt
\end{array}
\right.
\qquad \mbox{and} \qquad
\left\{
\begin{array}{llll}
Y_{k+\ell-1}=y^0_{k_{\ell}},\\
Y_{N+k+\ell-1}=\bar y^0_{k_{\ell}}.
\end{array}
\right.
\label{bY}
\end{eqnarray}
According to (\ref{M})-(\ref{bY}), the first order partial derivatives in (\ref{Form2_GradientJ_y0pq}) and (\ref{Form2_GradientJ_dot_y0pq}) are given by:
\begin{eqnarray}
\begin{array}{lll}
\hspace{-0.2cm} \bullet\; \displaystyle\frac{\partial J\big(Y\big)}{\partial Y_{k+\ell-1}}=\displaystyle\sum_{k=1}^K\displaystyle\sum_{\ell=1}^{m_k}\Big(M_{p+q-1,k + \ell -1}Y_{k+\ell-1} + M_{p+q-1,N + k + \ell -1}Y_{N + k+\ell-1}\Big) + b_{p+q-1},\\
\hspace{-0.2cm} \bullet \; \displaystyle\frac{\partial J\big(Y\big)}{\partial Y_{N+k+\ell-1}}=\displaystyle\sum_{k=1}^K\displaystyle\sum_{\ell=1}^{m_k}\Big(M_{N + p+q-1,k + \ell -1}Y_{k+\ell-1} + M_{N + p+q-1,N + k + \ell -1}Y_{N+k+\ell-1}\Big) + b_{N+p+q-1}.
\end{array}
\label{Lin_Syst_0}
\end{eqnarray}

Then one can write
\begin{eqnarray}
\nabla J\big(Y\big)=M Y + b.
\end{eqnarray}

\end{document}